%% file: arxiv.tex
\documentclass[journal]{IEEEtran}
\usepackage{graphics} 
\usepackage{epsfig} 
\usepackage{times} 
\usepackage{amsmath} 
\usepackage{amssymb}  
\usepackage{cite}
\usepackage{epstopdf}
\usepackage{subfigure,color,balance}
\usepackage{verbatim}
\usepackage{cases}
\usepackage{enumerate}

\usepackage{balance}

\usepackage[colorlinks=true,      
linkcolor=black,      
citecolor=black,      
filecolor=black,      
urlcolor=blue]{hyperref}

\newtheorem{definition}{Definition}
\newtheorem{theorem}{Theorem}
\newtheorem{proposition}{Proposition}
\newtheorem{remark}{Remark}
\newtheorem{lemma}{Lemma}
\newtheorem{corollary}{Corollary}
\newtheorem{problem}{Problem}
\newtheorem{step}{Step}

\input{defs}

\begin{document}
%
\title{
Controllability Analysis and Optimal Control of Mixed Traffic Flow with Human-driven and Autonomous Vehicles}

\author{Jiawei Wang, Yang Zheng, Qing Xu, Jianqiang Wang, and Keqiang Li

\thanks{This work is supported by National Key R\&D Program of China with 2016YFB0100906. The authors also acknowledge the support from TOYOTA. All correspondence should be sent to Y.~Zheng and K. Li.}
\thanks{J.W.~Wang, Q.~Xu, J.~Wang and K.~Li are with the School of Vehicle and Mobility, Tsinghua University, Beijing, China, and with the Center for Intelligent Connected Vehicles \& Transportation, Tsinghua University, Beijing, China. (wang-jw18@mails.tsinghua.edu.cn, \{qingxu,wjqlws,likq\}@tsinghua.edu.cn).}%
\thanks{Y. Zheng was with the Department of Engineering Science, University of Oxford, U.K. He is now with SEAS and CGBC at Harvard University ({zhengy@g.harvard.edu}).}%
}

\maketitle

\begin{abstract}
Connected and automated vehicles (CAVs) have a great potential to improve traffic efficiency in mixed traffic systems, which has been demonstrated by multiple numerical simulations and field experiments. However, some fundamental properties of mixed traffic flow, including controllability and stabilizability, have not been well understood. This paper analyzes the controllability of mixed traffic systems and designs a system-level optimal control strategy. Using the Popov-Belevitch-Hautus (PBH) criterion, we prove for the first time that a ring-road mixed traffic system with one CAV and multiple heterogeneous human-driven vehicles is not completely controllable, but is stabilizable under a very mild condition. Then, we formulate the design of a system-level control strategy for the CAV as a structured optimal control problem, where the CAV's communication ability is explicitly considered. Finally, we derive an upper bound for reachable traffic velocity via controlling the CAV. Extensive numerical experiments verify the effectiveness of our analytical results and the proposed control strategy. Our results validate the possibility of utilizing CAVs as mobile actuators to smooth traffic flow actively.

\end{abstract}


%
\IEEEpeerreviewmaketitle

\section{Introduction}
%
%
%
%

\IEEEPARstart{D}{uring}  the past decades, the increasing mobility demand has put a heavy burden on the existing transportation infrastructures. Designing advanced control methods has attracted significant research attention in order to improve traffic efficiency and road safety~\cite{baskar2011traffic}. Most existing strategies for traffic control rely on certain actuators at fixed locations, such as traffic signals and signs on roadside infrastructure~\cite{hegyi2005model}. Two typical systems are variable speed limits and variable speed advisory, which already have certain industrial applications~\cite{papageorgiou2008effects}. Due to their dependence on fixed infrastructure and drivers' compliance, however, the flexibility and effectiveness of these systems might be compromised~\cite{stern2018dissipation}.

As one key ingredient of traffic systems, the motion of vehicles plays an important role in traffic efficiency. Recent advancements on control and communication technologies have led to the emergence of connected and automated vehicles (CAVs), which are expected to revolutionize road transportation systems significantly. Compared to human-driven vehicles (HDVs), the cooperative formation of multiple CAVs, \eg, adaptive cruise control (ACC) and cooperative adaptive cruise control (CACC)~\cite{li2017dynamical}, has shown very promising effects on achieving higher traffic efficiency~\cite{talebpour2016influence}, better driving safety~\cite{fildes2015effectiveness} and lower fuel consumption~\cite{li2015effect}. These effects have been demonstrated and validated in some large-scale numerical simulations~\cite{schakel2010effects} and small-scale field experiments~\cite{milanes2014cooperative}. These technologies typically require that all the involved vehicles have autonomous capabilities. In practice, however, as the gradual deployment of CAVs, there will have to be a transition phase of mixed traffic systems containing both CAVs and HDVs. Due to the interactions between neighboring vehicles, it is possible to use a few CAVs as mobile actuators to influence the motion of their surrounding vehicles, which may in turn control the global traffic flow. This notion is known as the \emph{Lagrangian control} of traffic flow~\cite{stern2018dissipation,vinitsky2018lagrangian}. In this aspect, Cui \etal presented some theoretical analysis on the potential of one single CAV in a ring-road scenario~\cite{cui2017stabilizing}, and the pioneer real-world experiments in~\cite{stern2018dissipation} clearly demonstrated that one single CAV has a certain ability to dissipate stop-and-go waves in a mixed traffic flow. In this paper, we continue this direction of controlling traffic flow via CAVs, and present rigorous controllability analysis and optimal controller synthesis of mixed traffic systems consisting of one CAV and multiple heterogeneous HDVs.

Understanding the dynamics of mixed traffic systems is essential to reveal the full potential of CAVs. Many previous studies are based on numerical experiments for various traffic scenarios, which have demonstrated some positive effects of CAVs on improving traffic stability, capacity and throughput; see, \eg, \cite{shladover2012impacts,van2006impact}. These studies suggest that some improvement in traffic efficiency can be achieved if the penetration rate of CAVs reaches a certain level; otherwise, the existence of CAVs might lead to negligible influence~\cite{calvert2017will,kerner2016failure}. We note that these findings are valid only for specific control strategies and particular simulation setups. The theoretical potential of CAVs requires further investigation. Recently, a few results have been derived with respect to theoretical analysis of mixed traffic systems that reveal some fundamental properties. For example, Talebpour \etal \cite{talebpour2016influence} examined the string stability of a mixed vehicle platoon with infinite length and revealed a relationship between the traffic stability and the penetration rate of CAVs. The effect of CAVs' spatial distribution was investigated in~\cite{xie2018Heterogeneous}, which derived necessary conditions for linear stability of mixed traffic systems. More recently, theoretical controllability analysis was carried out by Cui \etal~\cite{cui2017stabilizing} and then extended by Zheng \etal~\cite{zheng2018smoothing}, where it is proved that mixed traffic flow in a ring-road scenario can be stabilized by controlling one single CAV. This stabilizability result reveals an essential ability of a single CAV in smoothing traffic flow. However, the theoretical results in~\cite{cui2017stabilizing,zheng2018smoothing} are only applicable to uniform traffic flow due to the homogeneous assumption for HDVs. In real traffic situations, it is necessary to consider the heterogeneous case where various types of drivers and vehicles coexist. This makes the controllability analysis more challenging, since the method of eigenvector calculation in~\cite{cui2017stabilizing} and the strategy of block diagonalization in~\cite{zheng2018smoothing} are not directly applicable.

In addition to analyzing the dynamics of mixed traffic systems, several new methods have been recently proposed for the control of CAVs in mixed traffic flow \cite{orosz2016connected,vaio2019Cooperative,jin2017optimal,stern2018dissipation,zheng2018smoothing,wu2018Stabilizing,wu2017flow,vinitsky2018lagrangian,guo2017an,pandita2013preceding,kesting2008adaptive,nishi2013theory,kreidieh2018dissipating}. One common feature of these methods is that the dynamics of other HDVs are considered in the system model explicitly. For example, a notion of connected cruise control (CCC), proposed by Orosz \etal~\cite{orosz2016connected}, exploits information from multiple HDVs ahead to design control decisions for the CAV at the tail of a vehicle string. In a similar setup, various topics have been investigated, including estimation of HDV dynamics~\cite{jin2018connected}, robustness against model uncertainties~\cite{jin2018experimental}, influence of control feedback gains~\cite{guo2017an}, and impairments of imperfect communication conditions~\cite{vaio2019Cooperative}. However, the role of communication abilities in mixed traffic flow has not been addressed explicitly in these studies. It is shown that different communication topologies have a big impact on the performance of a platoon formation; see, \eg,~\cite{zheng2016stability,zheng2018platooning}. In addition, the methods in~\cite{orosz2016connected,pandita2013preceding,vaio2019Cooperative,guo2017an} typically take local traffic performance around the CAV into account for controller design.

In this paper, we explicitly take into consideration the communication abilities of CAVs and highlight a crucial transformation towards the control goal of CAVs, from a local-level to a system-level. Specifically, one goal of CAV control is to achieve certain desired performance of the entire traffic system. This idea is indeed supported by previous research on dampening traffic waves. In this direction, the existing methods can be classified into three categories: 1) Heuristic methods \cite{stern2018dissipation,kesting2008adaptive,nishi2013theory}, which typically modify the CAV behavior based on different traffic states to dissipate traffic waves. For example, several extensions of ACC systems were developed to adapt their parameters based on current traffic situations~\cite{kesting2008adaptive}. A jam-absorption driving strategy based on geometric features of trajectory data was proposed in~\cite{nishi2013theory} to improve the effectiveness of CAV control. 
One disadvantage of this class of methods is that the parameter values may need to be tuned empirically in different traffic conditions, and the performance is not completely predictable. Moreover, they lack theoretical guarantees for their performance and tend to leave a long gap from its preceding vehicle, possibly causing other vehicles to cut in. 2) Learning-based methods~\cite{vinitsky2018lagrangian,wu2017flow,kreidieh2018dissipating}, which leverage machine learning frameworks, such as deep reinforcement learning, to train control strategies of CAVs. With improving the entire traffic flow as the training objective, the resulting strategy can enable CAVs to achieve a system-level control goal, for example, bottleneck decongestion~\cite{vinitsky2018lagrangian} and stop-and-go wave dissipation~\cite{kreidieh2018dissipating}. The shortages of these methods include that the non-transparent training process is usually computationally demanding, and the resulting strategies might lack generalizability and interpretability. 3) Model-based methods~\cite{zheng2018smoothing,wu2018Stabilizing,vaio2019Cooperative}, which adopt the perspective of rigorous control theory and offer certain insights for the CAV control problem in mixed traffic. 
For example, both Vaio \etal~\cite{vaio2019Cooperative} and Wu \etal~\cite{wu2018Stabilizing} formulated a controller synthesis problem that enables the CAV to dampen traffic waves based on string stability analysis. Recently, an optimal control method for CAVs was proposed in~\cite{zheng2018smoothing}, which attempts to achieve a system-level optimum criterion. Nevertheless, the requirement of global information of the entire traffic flow in~\cite{zheng2018smoothing} might restrict its practical applications.

In this paper, we focus on the controllability analysis and optimal controller synthesis of a mixed traffic system, which consists of one CAV and multiple heterogeneous HDVs in a ring road. First, we use a linearized car-following model to describe the behavior of mixed traffic flow. Then, the controllability analysis of this system is conducted based on the Popov-Belevitch-Hautus (PBH) criterion~\cite{skogestad2007multivariable}. Moreover, we formulate the design of a system-level control strategy for the CAV as a structured optimal control problem~\cite{jovanovic2016controller}, which incorporates explicit structural constraints according to the CAV's communication ability. Precisely, our contributions are as follows.
\begin{enumerate}
\item
We prove that a mixed traffic system with multiple heterogeneous HDVs and one single CAV is not completely controllable, but is stabilizable under a very mild condition. This result reveals a fundamental property of mixed traffic systems and confirms the feasibility of traffic control via CAVs, with no need of changing the behavior of HDVs. Also, our theoretical results validate the empirical observations in~\cite{stern2018dissipation} that a single CAV is able to stabilize mixed traffic flow and dampen undesired disturbances. In the case of homogeneous dynamics for HDVs, our result is consistent with~\cite{cui2017stabilizing,zheng2018smoothing}. Note that our proof relies heavily on the PBH test and eigenvalue-eigenvector analysis, and the methods in~\cite{cui2017stabilizing} and~\cite{zheng2018smoothing} are not directly suitable for the heterogeneous case.
\item
We formulate the problem of designing a control strategy for the CAV with limited communication abilities in the mixed traffic flow as structured optimal controller synthesis. Unlike existing methods, such as CACC~\cite{li2017dynamical} and CCC~\cite{orosz2016connected}, which attempt to achieve a better driving behavior of the CAV itself, the CAV control in our formulation directly aims to improve the performance of the entire traffic flow. In addition, the issue of limited communication is considered explicitly in our formulation. The structured optimal control is in general computationally intractable~\cite{jovanovic2016controller}, and we utilize the recent advance in sparsity invariance~\cite{furieri2019OnSep} to derive a convex relaxation, which allows one to compute a suboptimal solution efficiently.
\item
We present an analytical result with respect to the reachability property of the mixed traffic system, revealing a fundamental relationship between the desired state and the system dynamics. Based on the reachability analysis, we derive an upper bound of reachable velocity of mixed traffic flow and show how to design the desired system state in the controller. Extensive numerical experiments validate our theoretical results, and also confirm that the proposed controller can stabilize traffic flow and dampen traffic waves.

\end{enumerate}

\begin{figure}[t]
	\centering
	\subfigure[]
	{ \label{Fig:SystemModela}
		\includegraphics[scale=0.33, resolution=30]{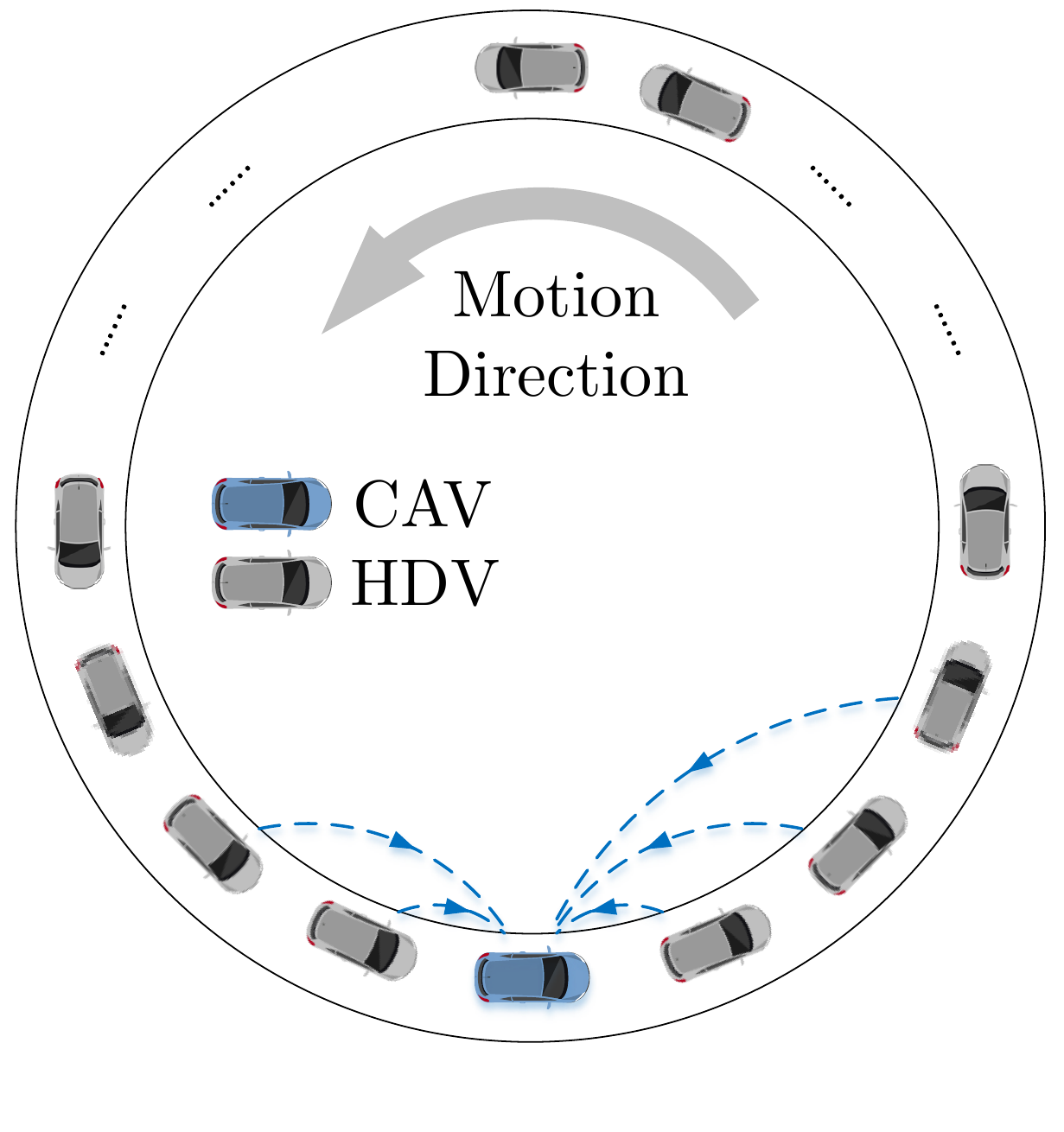}
	}
	\subfigure[]
	{\label{Fig:SystemModelb}
		\includegraphics[scale=0.33, resolution=30]{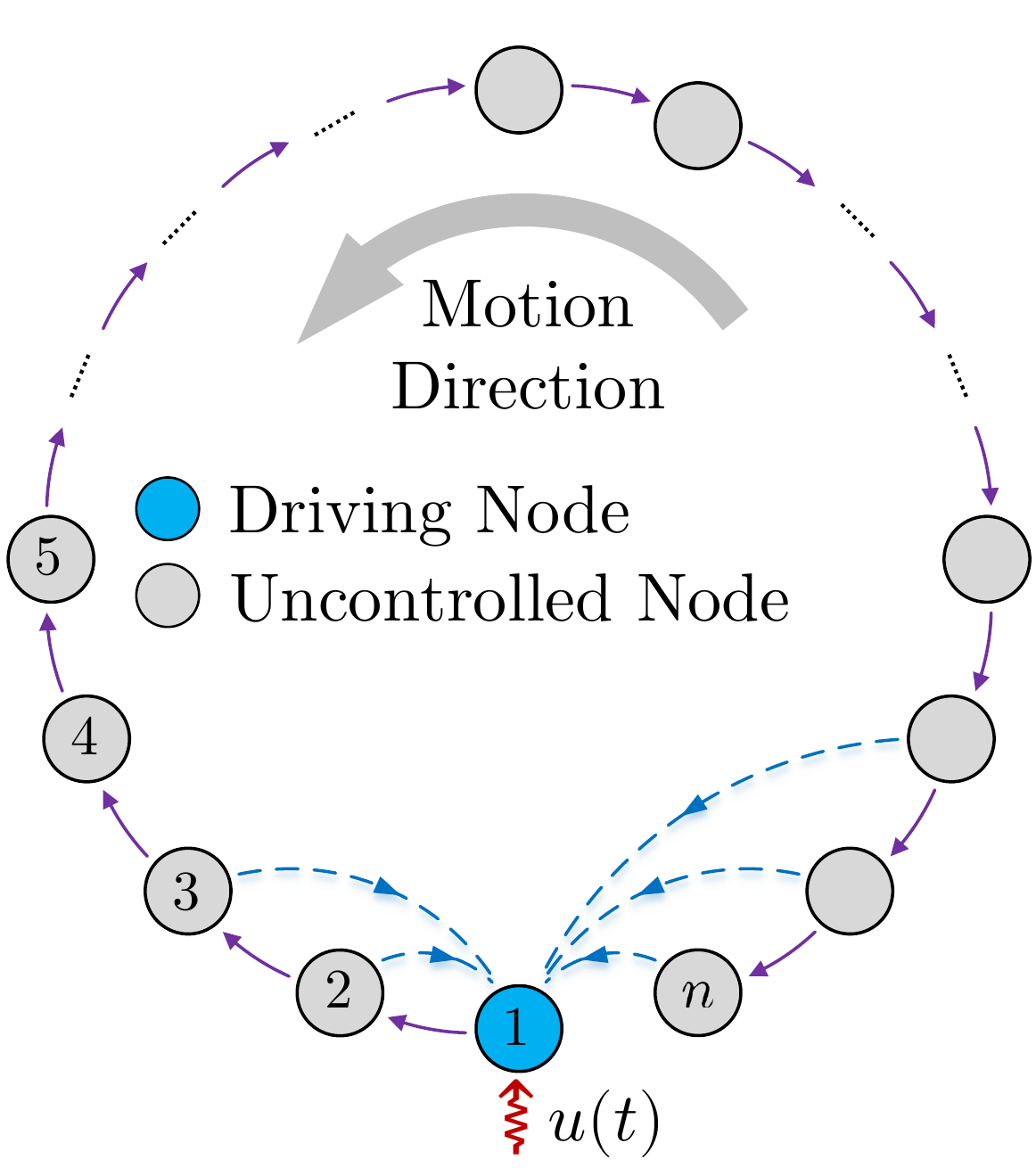}
	}
	\vspace{-2mm}
	\caption{System modeling schematic. (a) A single-lane ring road with one CAV and $n-1$ HDVs. (b) A simplified network system schematic. The blue arrows in (a)(b) represent the communication topology for the CAV; the purple arrows in (b) represent the interaction direction of car-following dynamics of HDVs.}
	\label{Fig:SystemModel}
\end{figure}

Some preliminary results have been presented in a conference version~\cite{wang2019controllability}. The rest of this paper is organized as follows. Section \ref{Sec:2}¢ò introduces the modeling for a mixed traffic system and the problem statement. The controllability result is presented in Section \ref{Sec:3}, and Section \ref{Sec:4} describes the theoretical framework to obtain a system-level optimal controller and the design of desired system state. Numerical simulations are shown in Section \ref{Sec:5}, and Section \ref{Sec:6} concludes this paper.

\section{System Modeling and Problem Statement}
\label{Sec:2}


\subsection{Modeling the Mixed Traffic System}
As shown in Fig.~\ref{Fig:SystemModel}, we consider a single-lane ring road with circumference $L$ consisting of one CAV and $n-1$ HDVs. It has been shown experimentally in \cite{sugiyama2008traffic,tadaki2013phase} that the ring road setup can easily reproduce the phenomenon of traffic waves without any infrastructure bottlenecks or lane changing behaviors. Moreover, as discussed in \cite{stern2018dissipation,cui2017stabilizing,zheng2018smoothing}, this setup also has several theoretical advantages, including 1) representative of a closed traffic system with no need for boundary conditions, 2) perfect control of average traffic density, 3) correspondence to an infinite straight road with periodic traffic dynamics.

\begin{figure}[t]
	\vspace{2mm}
	\centering
	\subfigure[]
	{\label{Fig:HDVModela}
		\includegraphics[scale=0.3]{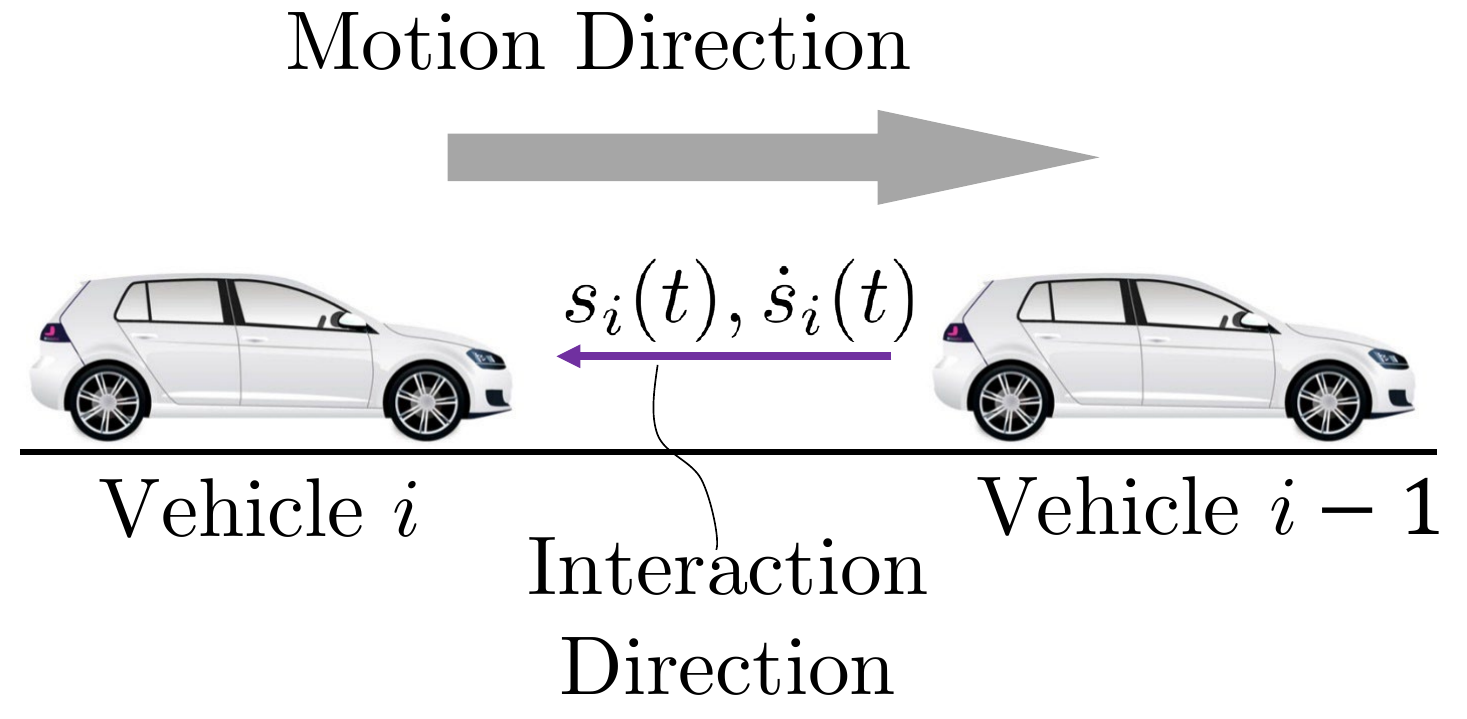}}
	\subfigure[]
	{\label{Fig:HDVRelationship}
		\includegraphics[scale=0.39]{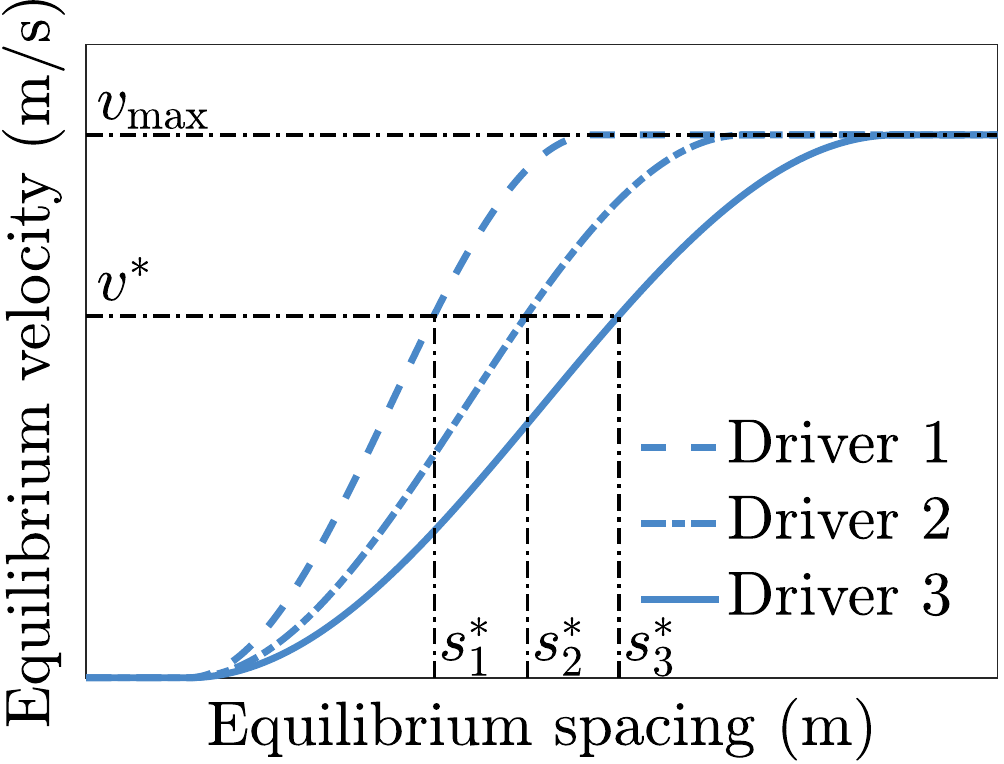}}
	\caption{Car-following dynamics for HDVs. (a) The driver usually considers the state of the preceding and his own vehicle. (b) A typical relationship between equilibrium spacing and equilibrium velocity, given by \eqref{Eq:Equilibrium}. Typically, $v^*$ grows as $s^*$ increases. For the same $v^*$, the equilibrium spacing $s^*$ for different drivers may be different.}
	\label{Fig:HDVModel}
\end{figure}

We denote the position, velocity and acceleration of vehicle $i$ as $p_i$, $v_i$ and $a_i$ respectively. The spacing of vehicle $i$, \ie, its relative distance from vehicle $i-1$, is defined as $s_i=p_{i-1}-p_i$. Without loss of generality, the vehicle length is ignored and we assume that vehicle no.1 is the CAV.
The optimal velocity model (OVM) \cite{bando1995dynamical} and intelligent driver model (IDM) \cite{Treiber2000Congested} are two typical models to describe car-following dynamics of human-driven vehicles. Both of them can be expressed as \cite{jin2017optimal}
\begin{equation}\label{Eq:HDVModel}
\dot{v}_i(t)=F_i\left(s_i(t),\dot{s}_i(t),v_i(t)\right),
\end{equation}
where $\dot{s}_i(t)=v_{i-1}(t)-v_i(t)$, and $F_i(\cdot)$ denotes that the acceleration of vehicle $i$ is determined by the relative distance, relative velocity and its own velocity; see Fig.~\ref{Fig:HDVModela} for illustration. In equilibrium traffic state, each vehicle moves with the same equilibrium velocity $v^*$, \ie, $v_i(t)=v^*, \dot{s}_i(t)=0$, for $i=1,\ldots,n$. Meanwhile, each vehicle has a corresponding equilibrium spacing $s_i^*$. According to \eqref{Eq:HDVModel}, $s_i^*$ of each HDV should satisfy
\begin{equation} \label{Eq:Equilibrium}
F_i\left(s_i^*,0,v^*\right)=0,\;i=2,\ldots,n.
\end{equation}
From \eqref{Eq:Equilibrium}, it is immediate to see that different choices of $v^*$ yield different values of $s_i^*$ for each HDV; see Fig.~\ref{Fig:HDVRelationship} for a typical relationship between $s_i^*$ and $v^*$. Unlike HDVs, the equilibrium spacing of the CAV, \ie, $s_1^*$, can be designed separately. The choice of $s_1^*$ is discussed in Section \ref{Sec:4c}.

Assuming that each vehicle has a small perturbation from the equilibrium state $(s_i^*,v^*)$, we define the error state between actual and equilibrium state of vehicle $i$ as
$$
x_i(t)=\begin{bmatrix}\tilde{s}_i(t),\tilde{v}_i(t)\end{bmatrix}^\tr=\begin{bmatrix}s_i(t)-s_i^*,v_i(t)-v^*\end{bmatrix}^\tr.
$$
Using \eqref{Eq:Equilibrium} and applying the first-order Taylor expansion to \eqref{Eq:HDVModel}, we can derive a linearized model for each HDV ($i=2,\ldots,n$)
\begin{equation}\label{Eq:LinearHDVModel}
\begin{cases}
\dot{\tilde{s}}_i(t)=\tilde{v}_{i-1}(t)-\tilde{v}_i(t),\\
\dot{\tilde{v}}_i(t)=\alpha_{i1}\tilde{s}_i(t)-\alpha_{i2}\tilde{v}_i(t)+\alpha_{i3}\tilde{v}_{i-1}(t),\\
\end{cases}
\end{equation}
with $\alpha_{i1} = \frac{\partial F_i}{\partial s_i}, \alpha_{i2} = \frac{\partial F_i}{\partial \dot{s}_i} - \frac{\partial F_i}{\partial v_i}, \alpha_{i3} = \frac{\partial F_i}{\partial \dot{s}_i}$ evaluated at the equilibrium state. To reflect the real driving behavior, we have $\alpha_{i1}>0, \alpha_{i2}>\alpha_{i3}>0$ \cite{cui2017stabilizing}. For the CAV, indexed as $i=1$, the acceleration signal is directly used as the control input $u(t)$, and its car-following model is
\begin{equation}\label{Eq:LinearCAVModel}
\begin{cases}
\dot{\tilde{s}}_1(t)=\tilde{v}_{n}(t)-\tilde{v}_1(t),\\
\dot{\tilde{v}}_1(t)=u(t).\\
\end{cases}
\end{equation}

To derive the global dynamics of the mixed traffic system, we lump the error states of all the vehicles into one global state, \ie, $x(t)=\begin{bmatrix}x_1^\tr (t),x_2^\tr (t),\ldots,x_n^\tr (t)\end{bmatrix}^\tr$. Then, a linearized state-space model for the mixed traffic system is obtained
\begin{equation} \label{Eq:SystemModel}
\dot{x}(t)=Ax(t)+Bu(t),
\end{equation}
where
$$
A=\begin{bmatrix} C_1 & 0 & \ldots &\ldots & 0 & C_2 \\
A_{22} & A_{21} & 0 & \ldots & \ldots & 0 \\
0 & A_{32} & A_{31} & 0 & \ldots & 0\\
\vdots & \ddots & \ddots & \ddots & \ddots & \vdots\\
0 & \ldots & \ldots & 0 & A_{n2} & A_{n1}\end{bmatrix},\;B=\begin{bmatrix} B_1 \\ B_2 \\B_2 \\ \vdots \\ B_2 \end{bmatrix},
$$
with each block matrix defined as
\begin{equation*}
\begin{aligned}
&A_{i1} = \begin{bmatrix} 0 & -1 \\ \alpha_{i1} & -\alpha_{i2} \end{bmatrix},
A_{i2} = \begin{bmatrix} 0 & 1 \\ 0 & \alpha_{i3} \end{bmatrix},\;i=2,3,\ldots,n\\
&C_1 = \begin{bmatrix} 0 & -1 \\ 0 & 0 \end{bmatrix},
C_2 = \begin{bmatrix} 0 & 1 \\ 0 & 0 \end{bmatrix},
B_1 = \begin{bmatrix} 0  \\ 1 \end{bmatrix},
B_2 = \begin{bmatrix} 0 \\ 0 \end{bmatrix}.
\end{aligned}
\end{equation*}

\begin{remark}
Note that unlike \cite{cui2017stabilizing,zheng2018smoothing,jin2017optimal} which focused on homogeneous dynamics only, we allow HDVs to have heterogeneous car-following dynamics $F_i(\cdot)$. Thus, the equilibrium spacing $s_i^*$ and the blocks $A_{i1},A_{i2}$ are in general different for different vehicles at the same equilibrium velocity $v^*$. This heterogeneity consideration is more suitable for practical scenarios, but also brings more challenges for theoretical analysis.
\end{remark}

\subsection{Problem Statement}

Before proceeding with the problem of designing an optimal control input $u(t)$ for the mixed traffic system \eqref{Eq:SystemModel}, we introduce two standard notions.

\begin{definition}[Controllability \rm{\cite{skogestad2007multivariable}}]
	The dynamical system $\dot{x}=Ax+Bu$, or the pair $(A,B)$, is controllable if, for any initial state $x(0)=x_0$, any time $t_f>0$ and any final state $x_f$, there exists an input $u(t)$ such that $x(t_f)=x_f$.
\end{definition}

It is well-known that if a system is controllable, we can place its closed-loop poles arbitrarily. For uncontrollable systems, there exist some uncontrollable modes that cannot be moved. If the uncontrollable modes are stable, we can still design a feedback controller to stabilize the closed-loop system. This leads to the notion of stabilizability.

\begin{definition}[Stabilizability \rm{\cite{skogestad2007multivariable}}]
	A system is stabilizable if its uncontrollable modes are all stable.
\end{definition}

Controllability and stabilizability are two fundamental properties of linear systems, which guarantee the existence of a stabilizing feedback controller. We note that many previous studies in mixed traffic systems only considered a specific CAV controller and focused on traffic stability, but have not discussed these notions explicitly. Two notable exceptions are in~\cite{cui2017stabilizing,zheng2018smoothing} where controllability analysis was conducted under the assumption of a homogeneous HDV model. In this paper, our first objective is to address the controllability and stabilizability of the heterogeneous mixed traffic system~\eqref{Eq:SystemModel}, formally stated as follows.
\begin{problem} [Controllability and Stabilizability]
	Investigate whether or under what circumstances the mixed traffic system \eqref{Eq:SystemModel} is controllable or stabilizable.
\end{problem}

\begin{figure}[t]
	\centering
	\includegraphics[scale=0.32]{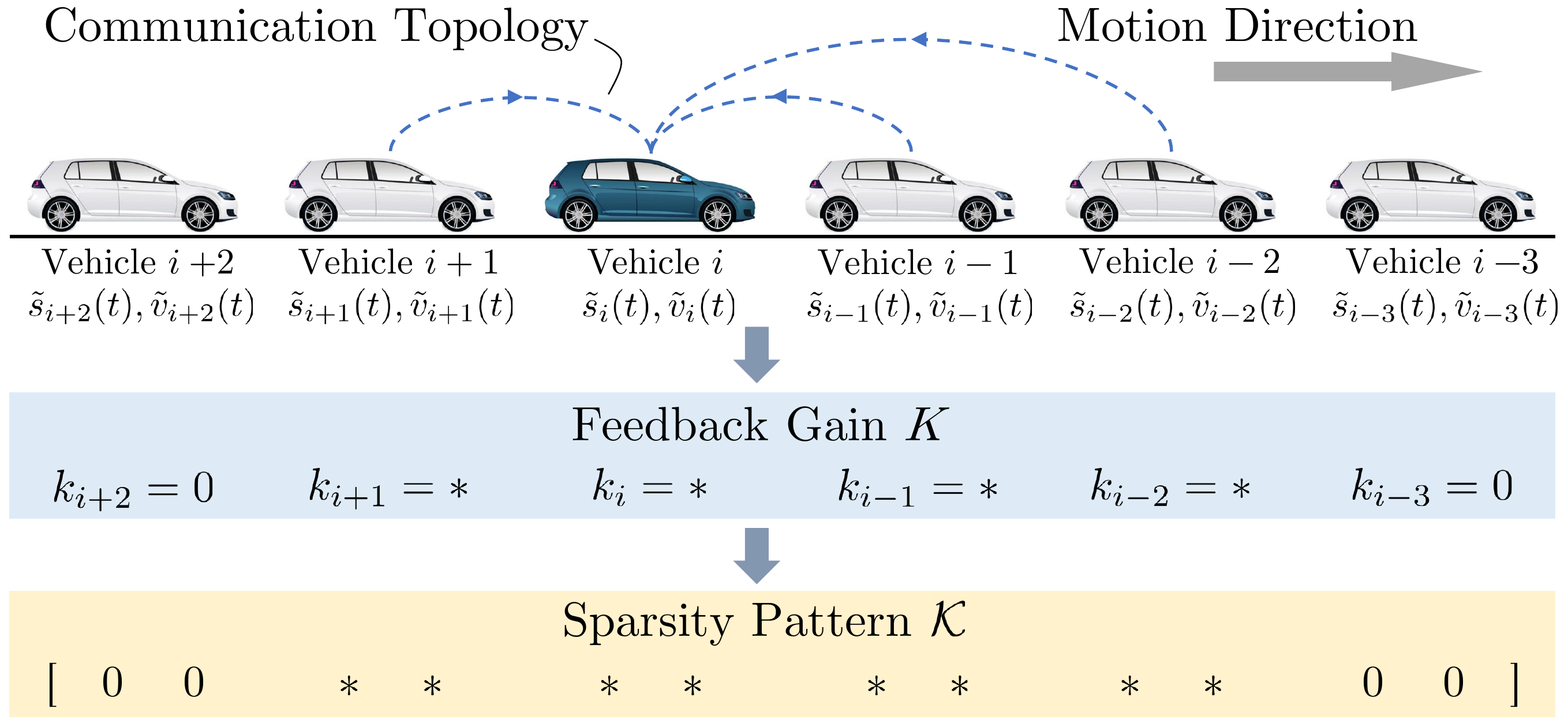}
	\caption{Illustration for structured constraints under limited communication abilities. The dashed blue arrows denote the communication topology of the CAV, \ie, vehicle $i$. Since the CAV can only have access to the state information of vehicles $i+1,i-1,i-2$, the specific feedback gain is zero for those vehicles from which vehicle $i$ cannot receive the information, \ie, vehicles $i+2,i-3$. The notation $*$ denotes that the specific feedback gain can be given a value of compatible dimensions.}
	\label{Fig:SparsityPattern}
\end{figure}

As depicted in Fig.~\ref{Fig:SystemModel}b, the mixed traffic flow can be viewed as a network system with the only CAV as a single driving node. Due to the limit of communication abilities in practice, the CAV can only receive partial information of the global traffic system for its feedback $u(t)$. Therefore, it is important to consider the local available information of the neighboring vehicles. Utilizing limited information exchange to control a large-scale network system leads to the notion of structured controller design~\cite{jovanovic2016controller}. To be precise, we define $\mathcal{E}^c\subseteq\{1,2,\ldots,n\}\times\{1\}$ as the communication network between the CAV and HDVs, where $(i,1)\in\mathcal{E}^c$ means that the CAV can receive the information from vehicle $i$. We consider a static state-feedback controller, \ie, $u(t)=-Kx(t)$, where $K=\begin{bmatrix}k_1,k_2,\ldots,k_n\end{bmatrix}\in\mathbb{R}^{1\times2n}$ and each block $k_i\in\mathbb{R}^{1\times2}$ represents the feedback gain of the state of vehicle $i$, \ie, $x_i(t)=\begin{bmatrix}\tilde{s}_i(t),\tilde{v}_i(t)\end{bmatrix}^\tr$. To reflect the communication topology, we require $k_i=0$, if $(i,1)\notin\mathcal{E}^c$. Hence, the communication topology requirement can be naturally imposed as a pre-specified sparsity pattern on $K$; see Fig.~\ref{Fig:SparsityPattern} for illustration. We define a block-sparsity pattern
$$
\mathcal{K}:=\{ K\in\mathbb{R}^{1\times2n}|k_i=0,\mathrm{if}\; (i,1)\notin\mathcal{E}^c,k_i\in\mathbb{R}^{1\times2} \}.
$$
Then, the second objective of this paper is as follows.
\begin{problem}[Optimal controller synthesis]
	compute a structured optimal controller $K\in\mathcal{K} $ for the CAV to dampen undesired perturbations in traffic flow, where $\mathcal{K}$ is determined by the communication topology $\mathcal{E}^c$.
\end{problem}

Perturbations in traffic flow are usually caused by lane changes and merges, but can also be generated in the absence of bottlenecks due to the stochastic nature of human's driving behavior~\cite{sugiyama2008traffic}. These undesired perturbations might result in traffic waves, which may easily cause series traffic congestion. In this paper, the undesired perturbation is modeled as a disturbance signal in each vehicle's acceleration. A precise formulation will be given in Section \ref{Sec:4a}.

\section{Controllability and Stabilizability}
\label{Sec:3}

In this section, we analyze the controllability and stabilizability of the mixed traffic system \eqref{Eq:SystemModel}.

\subsection{Controllability Analysis}

For completeness, we first introduce some notations. Denote $\dim (\cdot)$ as the dimension of a linear subspace, and $\ker M$ as the null space of $M \in \mathbb{R}^{n\times n}$, \ie, $\ker M=\{x\in \mathbb{R}^n| Mx=0\}$.

\begin{lemma}[Algebraic multiplicity \rm{\cite{horn2012matrix}}]
	\label{Lemma:Algebraic}
	Given an eigenvalue $\lambda$ of a matrix $M$, its algebraic multiplicity is the sum of the sizes of all its corresponding Jordan blocks. The number of its corresponding Jordan blocks is $\dim \, \ker (M-\lambda I)$ , and the number of its Jordan blocks of size $k$ is $2 \dim \, \ker (M-\lambda I)^k-\dim \, \ker (M-\lambda I)^{k+1}-\dim \, \ker (M-\lambda I)^{k-1}$.
\end{lemma}

\begin{lemma}[PBH controllability test \rm{\cite{skogestad2007multivariable}}]
	\label{Lemma:PBH}
	System $(A,B)$ is controllable, if and only if $\begin{bmatrix}\lambda I-A,B	\end{bmatrix}$ is of full row rank for all $\lambda$ being an eigenvalue of $A$.
\end{lemma}

It is clear that $\begin{bmatrix}\lambda I-A,B	\end{bmatrix}$ is of full row rank is equivalent to that the left null space of $\begin{bmatrix}\lambda I-A,B	\end{bmatrix}$ is empty, \ie, there exists no non-zero vector $\rho$ such that $\rho ^\tr \begin{bmatrix}\lambda I-A,B	\end{bmatrix}=0$. Conversely, if there exists $\rho \, (\rho \neq 0)$ such that $\rho ^\tr A=\lambda \rho ^\tr$ and $\rho ^\tr B=0$, then $(A,B)$ is not completely controllable, and $(\lambda,\rho)$ corresponds to an uncontrollable mode, which can also be given by $\rho^\tr x$ in the form of state component where $x$ is the state vector.

\begin{lemma}[Controllability invariance under state feedback \rm{\cite{skogestad2007multivariable}}]
	\label{Lemma:Invariance}
	$(A,B)$ is controllable, if and only if $(A-BK,B)$ is controllable for any state feedback $K\in \mathbb{R}^{1\times n}$ $(A\in \mathbb{R}^{n\times n})$. Furthermore, systems $(A,B)$ and $(A-BK,B)$ share the same uncontrollable modes.
\end{lemma}


\begin{figure}[t]
	\centering
	\includegraphics[scale=0.32]{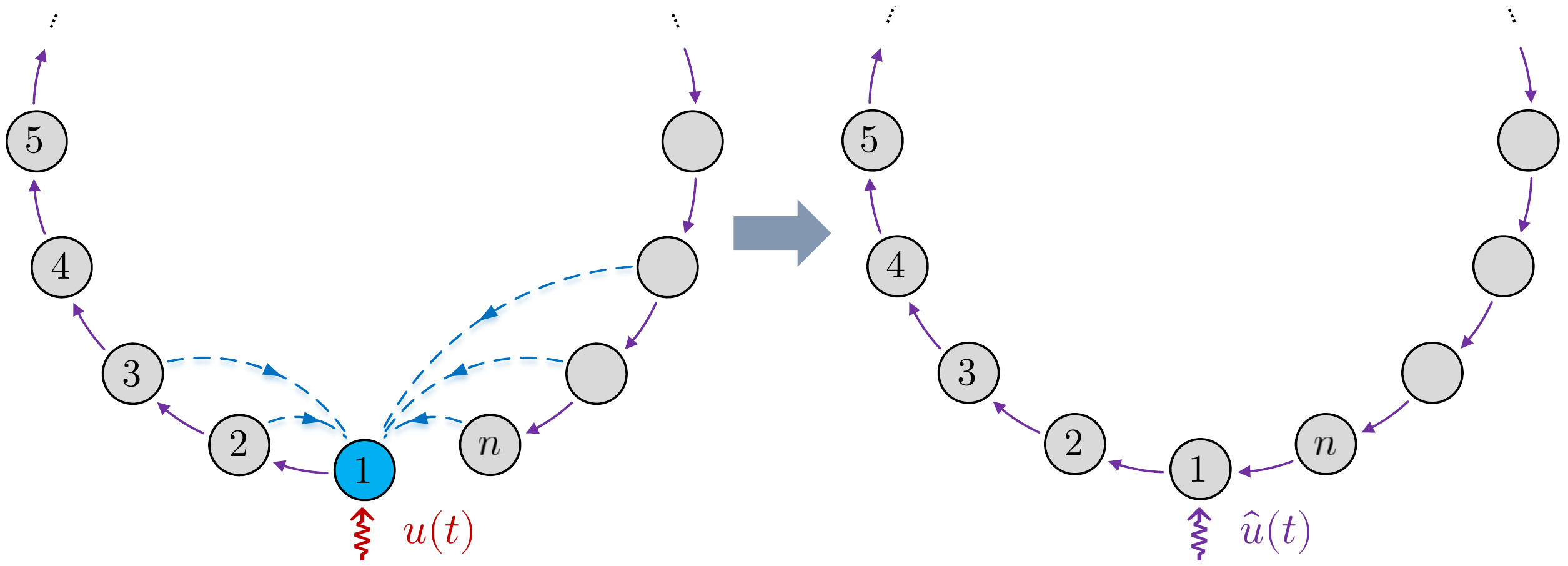}
	\caption{Illustration of system transformation and controllability invariance. Left and right correspond to $(A,B)$ in \eqref{Eq:SystemModel} and $(\hat A,B)$ in \eqref{Eq:SystemModelAfterTransformation} respectively. From left to right the CAV under control input $u(t)$ is replaced by an HDV with external control input $\hat u(t)$.}
	\label{Fig:Transformation}
\end{figure}

To facilitate the controllability analysis of the mixed traffic system \eqref{Eq:SystemModel}, we use Lemma \ref{Lemma:Invariance} and make the following transformation: suppose that the CAV has an external input
\begin{equation}
	\begin{aligned}
	\hat{u}(t) &= u(t)-\hat{K}x\\
	&=u(t)-\left(\alpha_{11}\tilde{s}_1(t)-\alpha_{12}\tilde{v}_1(t)+\alpha_{13}\tilde{v}_{n}(t)\right),
	\end{aligned}
\end{equation}
where $\hat{K}=\begin{bmatrix}\alpha_{11},-\alpha_{12},0,0,\ldots,0,\alpha_{13}\end{bmatrix}\in \mathbb{R}^{1\times 2n}$ and $\alpha_{11},\alpha_{12},\alpha_{13}$ are three positive constants. Then the original system $(A,B)$ can be transformed into $(\hat{A},B)$, described by the following dynamics
\begin{equation}
\label{Eq:SystemModelAfterTransformation}
\dot{x}(t)=\hat{A}x(t)+B\hat{u}(t),
\end{equation}
with $\hat{A}=A-B\hat{K}$, which can also be written as
\begin{equation}
	\label{Eq:Ahat}
\hat{A}=\begin{bmatrix} A_{11} & 0 & \ldots &\ldots & 0 & A_{12} \\
A_{22} & A_{21} & 0 & \ldots & \ldots & 0 \\
0 & A_{32} & A_{31} & 0 & \ldots & 0\\
\vdots & \ddots & \ddots & \ddots & \ddots & \vdots\\
0 & \ldots & \ldots & 0 & A_{n2} & A_{n1}\end{bmatrix},
\end{equation}
where
$$
A_{11} = \begin{bmatrix} 0 & -1 \\ \alpha_{11} & -\alpha_{12} \end{bmatrix},
A_{12} = \begin{bmatrix} 0 & 1 \\ 0 & \alpha_{13} \end{bmatrix}.
$$
$(\hat{A},B)$ describes a ring-road traffic system where all the vehicles are HDVs and one of them has an external control input $\hat{u}(t)$; see Fig.~\ref{Fig:Transformation} for illustration of the system transformation. Since $\hat{A}=A-B\hat{K}$, we have the following proposition according to Lemma \ref{Lemma:Invariance}.
\begin{proposition}
The controllability, stabilizability and uncontrollable modes remain the same between $(\hat{A},B)$ and $(A,B)$.
\end{proposition}

Therefore, we focus on system $(\hat{A},B)$ in the following. The first main result is as follows.

\begin{theorem} [Controllability]
	\label{Theorem:Controllability}
	Consider the mixed traffic system in a ring road with one CAV and $n-1$ heterogeneous HDVs given by \eqref{Eq:SystemModel}. The following statements hold:
	\begin{enumerate}
		\item
		System \eqref{Eq:SystemModel} is not completely controllable.
		\item
		There exists one uncontrollable mode corresponding to a zero eigenvalue, and this uncontrollable mode is stable.
	\end{enumerate}
\end{theorem}
	\begin{IEEEproof}
			Denote $ \rho_0=\begin{bmatrix}1,0,1,0,\ldots,1,0\end{bmatrix}^\tr\in \mathbb{R}^{2n\times1}$. Then it is easy to verify that $\rho_0^\tr \hat{A}=0\cdot \rho_0^\tr$ and $\rho_0^\tr B=0$, which means $\hat{A}$ has a zero eigenvalue and $\rank \begin{bmatrix}
			0 \cdot I-\hat{A},B
			\end{bmatrix}<2n$. According to Lemma \ref{Lemma:PBH}, we know that the zero eigenvalue corresponds to an uncontrollable mode.
			
			Next, we prove that the algebraic multiplicity of this zero eigenvalue is one. 
			According to Lemma \ref{Lemma:Algebraic}, we need to show that there is only one Jordan block of size one corresponding to the zero eigenvalue. Indeed, we can prove the following facts
				\begin{align}
					&\dim \, \ker (\hat{A}-0\cdot I)=1; \label{Eq:Algebraa}\\
					&\dim \, \ker (\hat{A}-0\cdot I)^2=1.\label{Eq:Algebrab}
				\end{align}
		Then by Lemma \ref{Lemma:Algebraic}, there is only one Jordan block corresponding to $\lambda=0$ and its size is $2\dim\,\ker \hat{A}^1-\dim\,\ker \hat{A}^2-\dim\,\ker \hat{A}^0=1$.
		To prove \eqref{Eq:Algebraa} and \eqref{Eq:Algebrab}, we consider the solution of $\hat{A}p=0$. Denote $p=\begin{bmatrix}	p_1^\tr,p_2^\tr,\ldots,p_n^\tr\end{bmatrix}^\tr$ where $p_i=\begin{bmatrix}
		p_{i1},p_{i2}\end{bmatrix}^\tr \in \mathbb{R}^{2\times 1},\,i=1,2,\ldots,n$. Considering the expression of $\hat{A}$ in \eqref{Eq:Ahat}, we know that $\hat{A}p=0$ is equivalent to $A_{i1} p_i+A_{i2} p_{i-1}=0$, leading to
		\begin{subequations}\label{Eq:ApzeroSolution}
			\begin{numcases}{}
			-p_{i2}+p_{(i-1)2}=0, \label{Eq:ApzeroSolutiona}\\
			\alpha_{i1}p_{i1}-\alpha_{i2}p_{i2}+\alpha_{i3}p_{(i-1)2}=0,\label{Eq:ApzeroSolutionb}
			\end{numcases}
		\end{subequations}
		for $i=1,2,\ldots,n$. For simplicity, we denote $i=n$ when $i$ takes the value of 0 and vice versa due to the circulant property. From \eqref{Eq:ApzeroSolutiona}, we know $p_{i2}=c$, $i=1,\ldots,n$ where $c$ is a constant. Substituting it into \eqref{Eq:ApzeroSolutionb} leads to $p_{i1}=\frac{\alpha_{i2}-\alpha_{i3}}{\alpha_{i1}}  c$. Then, the solution of $\hat{A}p=0$ is of the unique form
		\begin{equation}
		\label{Eq:Solutionofp}
			p=c\begin{bmatrix}
			\frac{\alpha_{12}-\alpha_{13}}{\alpha_{11}},1,\ldots,\frac{\alpha_{n2}-\alpha_{n3}}{\alpha_{n1}},1
			\end{bmatrix}^\tr,
		\end{equation}
		which has dimension one. Therefore, the statement in \eqref{Eq:Algebraa} holds. The proof of \eqref{Eq:Algebrab} is similar, and details can be found in Appendix \ref{Sec:AppendixA}.
		
	Now, we know that the algebraic multiplicity of the zero eigenvalue is one, and $\rho_0$ is the only mode corresponding to the zero eigenvalue. Accordingly, the uncontrollable mode $\rho_0$ will remain constant, \ie, equal to its initial state, in the system evolution. Therefore, we can conclude that the uncontrollable mode is stable in the Lyapunov sense.

	\end{IEEEproof}

 The uncontrollable mode can be expressed as $\rho _{0}^\tr x \left( t \right) = \sum _{i=1}^{n}\tilde{s}_{i} \left( t \right) = \sum _{i=1}^{n}s_{i} \left( t \right) - \sum _{i=1}^{n}s_{i}^{*}$. 
The physical interpretation is that the uncontrollable mode reflects the ring-road setup in system \eqref{Eq:SystemModel}, indicating that the sum of each vehicle's spacing must remain constant, \ie, $\sum _{i=1}^{n}s_{i} \left( t \right) =L $. Recall that the equilibrium spacing of the CAV, \ie, $s_1^*$, can be designed separately. Therefore, if we have $s_{1}^{*}=L- \sum _{i=2}^{n}s_{i}^{*} $, where $s_{2}^{*} , \ldots ,  s_{n}^{*} $ are given by \eqref{Eq:Equilibrium}, then $\rho_0^\tr x(t)=0$ for all $t\geq 0$. The influence of different choices of $s_1^*$ is further discussed in Section \ref{Sec:4c}.

\subsection{Stabilizability Analysis}
	Here, we show that all other modes corresponding to non-zero eigenvalues of $\hat{A}$ are controllable under a mild condition, which leads to the following theorem.
	
	\begin{theorem}[Stabilizability]
		\label{Theorem:Stabilizability}
		Consider the mixed traffic system in a ring road with one CAV and $n-1$ heterogeneous HDVs given by \eqref{Eq:SystemModel}. System \eqref{Eq:SystemModel} is stabilizable, if the following condition holds,
			\begin{equation} \label{Eq:ControllabilityTheorem}
		\alpha_{j1}^2-\alpha_{i2}\alpha_{j1}\alpha_{j3}+\alpha_{i1}\alpha_{j3}^2 \neq 0,\; \forall i,j\in\{1,2,\ldots,n\}.
		\end{equation}
	\end{theorem}
	
	\begin{IEEEproof}
		We first show that condition \eqref{Eq:ControllabilityTheorem} yields the following fact: each non-zero eigenvalue $\lambda$ of $\hat{A}$ satisfies
		\begin{equation} \label{Eq:Eigenvalue}
		\lambda ^{2}+ \alpha _{i2} \lambda + \alpha _{i1} \neq 0,\;   \alpha _{i3} \lambda + \alpha _{i2} \neq 0,\; \forall i \in  \{ 1,2, \ldots ,n \},
		\end{equation}
		The proof is elementary and can be found in Appendix \ref{Sec:AppendixB}.
		
		Assume that there exists a non-zero eigenvalue $\lambda$ of $\hat{A}$ which corresponds to an uncontrollable mode. Then according to Lemma \ref{Lemma:PBH}, there exists a left eigenvector $\rho$ of $\hat{A}$ associated with $\lambda$ such that $\rho^\tr B=0$. Denote $\rho$ as
		\begin{equation} \label{Eq:FormofRho}
			\rho = \begin{bmatrix}
			\rho _{1}^\tr , \rho _{2}^\tr , \ldots , \rho _{n}^\tr
			\end{bmatrix}  ^\tr ,
		\end{equation}
		where $\rho _{i}= \left[  \rho _{i1}, \rho _{i2} \right] ^\tr  \in \mathbb{R}^{2\times 1} $. Considering that only the second element in $B$ is non-zero, this assumption leads to $ \rho _{12}=0 $. In addition, due to $  \rho ^\tr A= \lambda  \rho ^\tr  $, we know for $ i=1,2, \ldots ,n $,
		$$
		\rho _{i}^\tr  \left(  \lambda I-A_{i1} \right) = \rho _{i+1}^\tr A_{ \left( i+1 \right) 2}
		$$
		Since $  \lambda I-A_{i1} $ is invertible according to \eqref{Eq:Eigenvalue}, we then have
		\begin{equation}\label{Eq:Eigenvector}
		\rho _{i}^\tr = \rho _{i+1}^\tr A_{ \left( i+1 \right) 2} \left(  \lambda I-A_{i1} \right) ^{-1},
		\end{equation}
		for $ i=1,2, \ldots ,n $. Upon denoting $ D_{i1}= \left(  \lambda I-A_{i1} \right) ^{-1}$, and considering the recursive result of \eqref{Eq:Eigenvector}, we have
		\begin{equation}\label{Eq:Step2Cal}
		\begin{aligned}
		\rho _{1}^\tr &= \rho _{2}^\tr A_{22}D_{11} \\
		&= \rho _{3}^\tr A_{32}D_{21}A_{22}D_{11}  \\
		&= \cdots  \\
		&= \rho _{n}^\tr A_{n2}D_{ \left( n-1 \right) 1}A_{ \left( n-1 \right) 2} \cdots D_{21}A_{22}D_{11}\\
		&= \rho _{1}^\tr A_{12} \left( D_{n1}A_{n2} \right)  \cdots  \left( D_{21}A_{22} \right) D_{11},
		\end{aligned}
		\end{equation}
		where $ D_{i1}A_{i2}$ is equal to 
		\begin{equation}\label{Eq:Step2Value}
		D_{i1}A_{i2}=\frac{1}{a_{i}} \begin{bmatrix}
		0  &  c_{i}\\
		0  &  b_{i}\\
		\end{bmatrix}, i=1,2, \ldots ,n,
		\end{equation}
		with $ a_{i}= \lambda ^{2}+ \alpha _{i2} \lambda + \alpha _{i1}$, $ b_{i}= \alpha _{i3} \lambda + \alpha _{i2}$, and $ c_{i}= \lambda + \alpha _{i2}- \alpha _{i3} $.
		According to \eqref{Eq:Eigenvalue}, we have $ a_{i} \neq 0$ and $ b_{i} \neq 0 $. Substituting \eqref{Eq:Step2Value} into \eqref{Eq:Step2Cal} yields
		\begin{equation}
		\begin{bmatrix}
		\rho _{11}  &   \rho _{12}
		\end{bmatrix}
		=  \begin{bmatrix}
		\rho _{11}  &   \rho _{12}
		\end{bmatrix}
		\begin{bmatrix}
		\alpha _{11}  &   \lambda \\
		\alpha _{11} \alpha _{13}  &   \lambda  \alpha _{13}\\
		\end{bmatrix}
		\frac{ \Pi _{i=2}^{n}b_{i}}{ \Pi _{i=1}^{n}a_{i}}.
		\end{equation}
		Since $\rho _{12}=0 $, we know $ \lambda  \rho _{11} \Pi _{i=2}^{n}b_{i}=0$. Because $ \lambda  \neq 0$ and $ b_{i} \neq 0$ , we have  $\rho _{11}=0 $, leading to $ \rho _{1}=0 $. Using \eqref{Eq:Eigenvector} recursively, we can obtain
		$
		\rho = \begin{bmatrix}
		\rho _{1}^\tr , \rho _{2}^\tr , \ldots , \rho _{n}^\tr
		\end{bmatrix}^\tr =0,
		$
		which contradicts $\rho\neq 0$. Consequently, the assumption does not hold. Therefore, all the modes corresponding to non-zero eigenvalues are controllable.
		
		Recall that we have shown in Theorem \ref{Theorem:Controllability} that there is only one uncontrollable mode $\rho_0^\tr x(t)$, which corresponds to the zero eigenvalue and is stable. Therefore, we conclude that the mixed traffic system \eqref{Eq:SystemModel} is stabilizable under condition \eqref{Eq:ControllabilityTheorem}.
	\end{IEEEproof}

	Note that \eqref{Eq:ControllabilityTheorem} is a sufficient condition for the mixed traffic system to be stabilizable. This condition restricts the locations of the closed-loop poles as shown in \eqref{Eq:Eigenvalue}, and allows us to reveal the characteristic of their corresponding eigenvectors. When choosing $\alpha_{i1}$, $\alpha_{i2}$, $\alpha_{i3}$ randomly, condition \eqref{Eq:ControllabilityTheorem} is satisfied with probability one; accordingly, the heterogeneous mixed traffic system \eqref{Eq:SystemModel} is stabilizable with probability one.
	
	In \cite{zheng2018smoothing}, Kalman controllability criterion was leveraged to analyze the controllability property of mixed traffic systems with homogeneous HDVs. However, for the heterogeneous mixed traffic system \eqref{Eq:SystemModel}, it is nontrivial to compute the rank of the Kalman controllability matrix $Q_c = \begin{bmatrix}
	B,AB,\ldots,A^{2n-1}B	\end{bmatrix}$. Instead, the eigenvalue-eigenvector analysis based on the PBH test gives analytical solutions for the controllability and stabilizability of the heterogeneous case.	In the case of homogeneous traffic flow, Theorems \ref{Theorem:Controllability} and~\ref{Theorem:Stabilizability} are consistent with that in \cite{zheng2018smoothing}. 
Specifically, if~$\alpha_{i1}=\alpha_1,\alpha_{i2}=\alpha_2,\alpha_{i3}=\alpha_3$, for $i=1,2,\ldots,n$, condition \eqref{Eq:ControllabilityTheorem} reduces to $\alpha_1-\alpha_2 \alpha_3+\alpha_3^2 \neq 0$. As proved in \cite{zheng2018smoothing}, if $\alpha_1-\alpha_2 \alpha_3+\alpha_3^2 = 0$, the homogeneous mixed traffic system is still stabilizable since all the uncontrollable modes are asymptotically stable.
	
	\begin{remark}
		Some previous works imposed a lower bound on the CAV penetration rate to guarantee the stability of mixed traffic flow; see \eg, \cite{talebpour2016influence,van2006impact,xie2018Heterogeneous,wu2018Stabilizing}. We note that these works usually employed specific controllers for CAVs, such as ACC- or CACC-type controllers, which might restrict the potential of CAVs. Instead, if one focuses on stabilizability of mixed traffic flow directly, Theorem \ref{Theorem:Stabilizability} has no requirement on system size $n$ or the stability property of the original traffic system with HDVs only. We prove that the entire traffic flow can be stabilized by controlling one single CAV, showing a greater capability of CAVs beyond normal expectations. Indeed, our theoretical results validate the empirical observations from field experiments in \cite{stern2018dissipation} and the training results from reinforcement learning in \cite{wu2017flow}, that one single CAV can be used as a mobile actuator to dampen traffic waves.
	\end{remark}

\section{Optimal Controller Synthesis}
\label{Sec:4}
The stabilizability of the mixed traffic system \eqref{Eq:SystemModel} guarantees the existence of a stabilizing controller. In this section, we proceed to design an optimal control input $u(t)$ for the CAV to dampen undesired perturbations of the traffic flow. We first formulate this task as a structured optimal control problem, and then introduce a numerical solution approach based on convex relaxation. We also discuss how to design the desired system state based on reachability analysis.

\subsection{System-level Performance and Structured Optimal Control}
\label{Sec:4a}
Previous works on the control of CAVs usually focused on local-level performance, \ie, to achieve a better driving behavior of a single CAV \cite{milanes2014cooperative} or a CAV platoon \cite{li2017dynamical}. Here, we first define a system-level performance index for the entire traffic system \eqref{Eq:SystemModel}. In particular, we aim to minimize the influence of undesired perturbations on traffic flow by controlling one CAV. To model this scenario, we assume that there exist certain disturbances
\begin{equation}
\label{Eq:Distrubance}
	\omega (t) = \begin{bmatrix}
	\omega_1(t),\omega_2(t),\ldots,\omega_n(t)
	\end{bmatrix}^\tr,
\end{equation}
where $\omega_i (t)$ is a scalar disturbance signal with finite energy in the acceleration of vehicle $i$ ($i=1,2,\ldots,n$). Then, the system model \eqref{Eq:SystemModel} becomes
\begin{equation}\label{Eq:DisturbanceModel}
\dot{x}(t)=Ax(t)+Bu(t)+H\omega(t),
\end{equation}
where $H\in\mathbb{R}^{2n\times n}$ is
$$
H=\begin{bmatrix}
H_{1}  &  0  &   \cdots   &  0\\
0  &  H_{1}  &  \ddots  &   \vdots \\
\vdots   &  \ddots  &  \ddots  &  0\\
0  &   \cdots   &  0  &  H_{1}\\
\end{bmatrix},
$$
with the block entry denoting $ H_{1}= \left[ 0,1 \right] ^{\tr} $. We use
\begin{equation}
\label{Eq:PerformanceOutput}
	z \left( t \right) =\begin{bmatrix}
	\gamma _{s}\tilde{s}_{1} \left( t \right) , \gamma _{v}\tilde{v}_{1} \left( t \right) , \ldots , \gamma _{s}\tilde{s}_{n} \left( t \right) , \gamma _{v}\tilde{v}_{n} \left( t \right) , \gamma _{u}u \left( t \right)
	\end{bmatrix}  ^\tr
\end{equation}
to denote a performance output, with weight coefficients $\gamma _{s}, \gamma _{v}, \gamma _{u}>0$ representing the penalty for spacing error, velocity error and control input, respectively. This system-level performance output can also be written as
\begin{equation}\label{Eq:Output}
z(t)=\begin{bmatrix}
Q^{\frac{1}{2}}\\0
\end{bmatrix}x(t)+\begin{bmatrix}
0\\R^{\frac{1}{2}}
\end{bmatrix}u(t),
\end{equation}
with  $Q^{\frac{1}{2}}=\mathrm{diag} \left(  \gamma _{s}, \gamma _{v}, \ldots , \gamma _{s}, \gamma _{v} \right)$, $R^{\frac{1}{2}}= \gamma _{u} $. In \eqref{Eq:Distrubance} and \eqref{Eq:PerformanceOutput}, we allow the perturbation to arise from anywhere in the traffic flow, and the performance output $z(t)$ takes into account all the vehicles' deviations in the traffic flow. This setup indicates a system-level consideration.

We use $G_{z \omega }$ to denote the transfer function from disturbance  $\omega$  to performance output $z $. Upon denoting $  \Vert G_{z \omega} \Vert $  as the  $\mathcal{H}_{2}$  norm that reflects the influence of disturbances, the design of an optimal control input $u \left( t \right)$ for the CAV under a pre-specified communication topology can be formulated as
\begin{equation}\label{Pr:Original}
\begin{aligned}
\min_{K} \quad &\Vert G_{z \omega} \Vert^2\\
\mathrm{subject\;to} \quad &u=-Kx,\;K\in \mathcal{K},
\end{aligned}
\end{equation}
where $K\in\mathcal{K}$ reflects the information that is available for the CAV. Note that previous local-level controllers, \eg, CACC \cite{li2017dynamical} or CCC \cite{orosz2016connected}, mostly only utilized the state information from vehicles ahead of the CAV to determine its own behavior. In this paper, our system-level consideration directly aims to improve the entire traffic system by controlling the CAV. Therefore, the CAV can utilize information from any vehicle within the communication range, \eg, the vehicles ahead and the vehicles behind.

Formulation \eqref{Pr:Original} is known as the structured optimal control problem \cite{jovanovic2016controller,zheng2017convex}. This problem is in general non-convex and computationally hard to find a globally optimal solution. In the following, we highlight that one particular difficulty lies in $K\in\mathcal{K}$.
\begin{lemma} [$\mathcal{H}_{2}$ norm of a transfer function \rm{\cite{skogestad2007multivariable}}]
	\label{Lemma:H2norm}
	Consider a stable system with the dynamics $\dot{ x} \left( t \right) =Ax \left( t \right) +H\omega \left( t \right)$, and the output $ z \left( t \right) =Cx \left( t \right)  $. The $ \mathcal{H}_{2} $ norm of the transfer function from $ w \left( t \right) $ to $ z \left( t \right) $ can be computed by
	\begin{equation}
	\Vert G_{z \omega} \Vert ^{2}=\inf_{X\succ0} \{ \mathrm{Tr} ( CXC^\tr )  \vert AX+XA^\tr +HH^\tr \preceq 0 \}.
	\end{equation}
\end{lemma}
\vspace{1mm}

Lemma \ref{Lemma:H2norm} offers a standard technique to compute the $\mathcal{H}_2$ norm of the transfer function of a linear system. Using $u=-Kx$, \eqref{Eq:DisturbanceModel} and \eqref{Eq:Output} become
\begin{equation} \label{Eq:CompactFormofxandz}
\begin{aligned}
&\dot{x}(t)=(A-BK)x(t)+H\omega(t),\\
&z(t)=\begin{bmatrix}
Q^{\frac{1}{2}}\\-R^{\frac{1}{2}}K
\end{bmatrix}x(t).
\end{aligned}
\end{equation}
Based on \eqref{Eq:CompactFormofxandz} and Lemma \ref{Lemma:H2norm}, problem \eqref{Pr:Original} can be equivalently reformulated as
\begin{equation}\label{Pr:P1}
\begin{aligned}
\min_{X,K} \quad &\mathrm{Tr}(QX)+\mathrm{Tr}(K^\tr RKX)\\
\mathrm{subject\;to}\quad & (A-BK)X+X(A-BK)^\tr+HH^\tr\preceq 0,\\
&X\succ 0,\;K\in\mathcal{K}.
\end{aligned}
\end{equation}
Using a standard change of variables \cite{jovanovic2016controller,zheng2017convex}
\begin{equation}
K=ZX^{-1},
\end{equation}
pre-and post-multiplying $AX+XA^\tr+HH^\tr \preceq 0$ by $P=X^{-1}$, and using the Schur complement, we can obtain an equivalent form of \eqref{Pr:P1}
\begin{equation}\label{Pr:P2}
\begin{aligned}
\min_{X,Y,Z} \quad &\mathrm{Tr}(QX)+\mathrm{Tr}(RY)\\
\mathrm{subject\;to}\quad & AX+XA^\tr-BZ-Z^\tr B^\tr+HH^\tr\preceq 0,\\
&\begin{bmatrix}
Y&Z\\Z^\tr&X
\end{bmatrix}\succeq 0,\;X\succ 0,\;ZX^{-1}\in\mathcal{K}.
\end{aligned}
\end{equation}
In the absence of $ZX^{-1}\in\mathcal{K}$, problem \eqref{Pr:P2} will be convex and can be solved via existing conic solvers, \eg, Mosek~\cite{andersen2000mosek}. However, the constraint $K\in\mathcal{K}$ appears naturally in the controller design for the CAV in a mixed traffic system, since the CAV can only use its neighboring vehicles' information for feedback.

	\begin{remark}
		In the proposed optimal strategy, the setup of the weight coefficients $\gamma _{s}, \gamma _{v}, \gamma _{u}$ play a significant role in the control performance. 
In principle, a larger value of $\gamma_{s}, \gamma_{v}$ indicates more consideration of the state error of all the vehicles in traffic flow, which normally allows the CAV to stabilize traffic flow in a shorter time. Conversely, a larger value of $\gamma _u$ commonly keeps a lower control input of the CAV, leading to more feasible control actions within practical bounds. Existing heuristic controllers, \eg, FollowerStopper and PI with Saturation \cite{stern2018dissipation}, have 
parameters that need to be adjusted empirically, and the resulting performance is not completely predictable. In contrast, our strategy shows an evident advantage in tuning the controller, which can achieve preferable performance or adapt to different traffic conditions. Note that we utilize a homogeneous setup for weight coefficients of all the vehicles in performance output \eqref{Eq:PerformanceOutput}. In fact, heterogeneity can also be introduced to show different levels of consideration for different vehicles, and Formulation \eqref{Pr:P2} still works to return an optimal feedback gain.
	\end{remark}

\begin{remark}
	Without the constraint $K\in\mathcal{K}$,  \eqref{Pr:Original} is a standard $\mathcal{H}_{2}$ optimal control problem, for which efficient methods are available to compute an optimal solution. Note that in \eqref{Pr:P2}, the solution $X$ can build a Lyapunov function  $V(x)=x^\tr (t) X^{-1} x(t)$ for the closed-loop system. The structured optimal control problem \eqref{Pr:Original} (or its variants) has attracted some attention in the literature. A few methods have been proposed to find an approximation solution, such as using convex approximations \cite{zheng2017convex}, or directly employing non-convex optimization techniques~\cite{jovanovic2016controller}. However, many existing methods require that the system is completely controllable, and therefore they are not applicable to our problem since the mixed traffic system is not completely controllable, as proved in Theorem~\ref{Theorem:Controllability}.
\end{remark}

\subsection{Numerical Solution Approach}
\label{Sec:4b}
In this section, we utilize a recent strategy based on sparsity invariance, originally introduced in \cite{furieri2019OnSep}, to compute a specific control feedback gain $K\in \mathcal{K}$ for problem \eqref{Pr:Original}.

To highlight the idea of sparsity invariance, we introduce a few notations to represent sparse matrices. Given a matrix $M\in \mathbb{R}^{m\times n}$, $M_{ij}$ denotes its entry in $i$-th row, $j$-th column. We denote $\{0,1\}^{m\times n}$ as the set of $m\times n$ binary matrices, meaning each entry in the matrix is either 0 or 1. Given $M\in \{0,1\}^{m\times n}$, $\mathrm{Sparse}(M)$ denotes a corresponding sparsity pattern, defined as
$$
\mathrm{Sparse}(M):=\{A\in \mathbb{R}^{m\times n}|A_{ij}=0,\,\mathrm{if}\, M_{ij}=0 \}.
$$

\begin{figure}[t]
	\centering
	\includegraphics[scale=0.22]{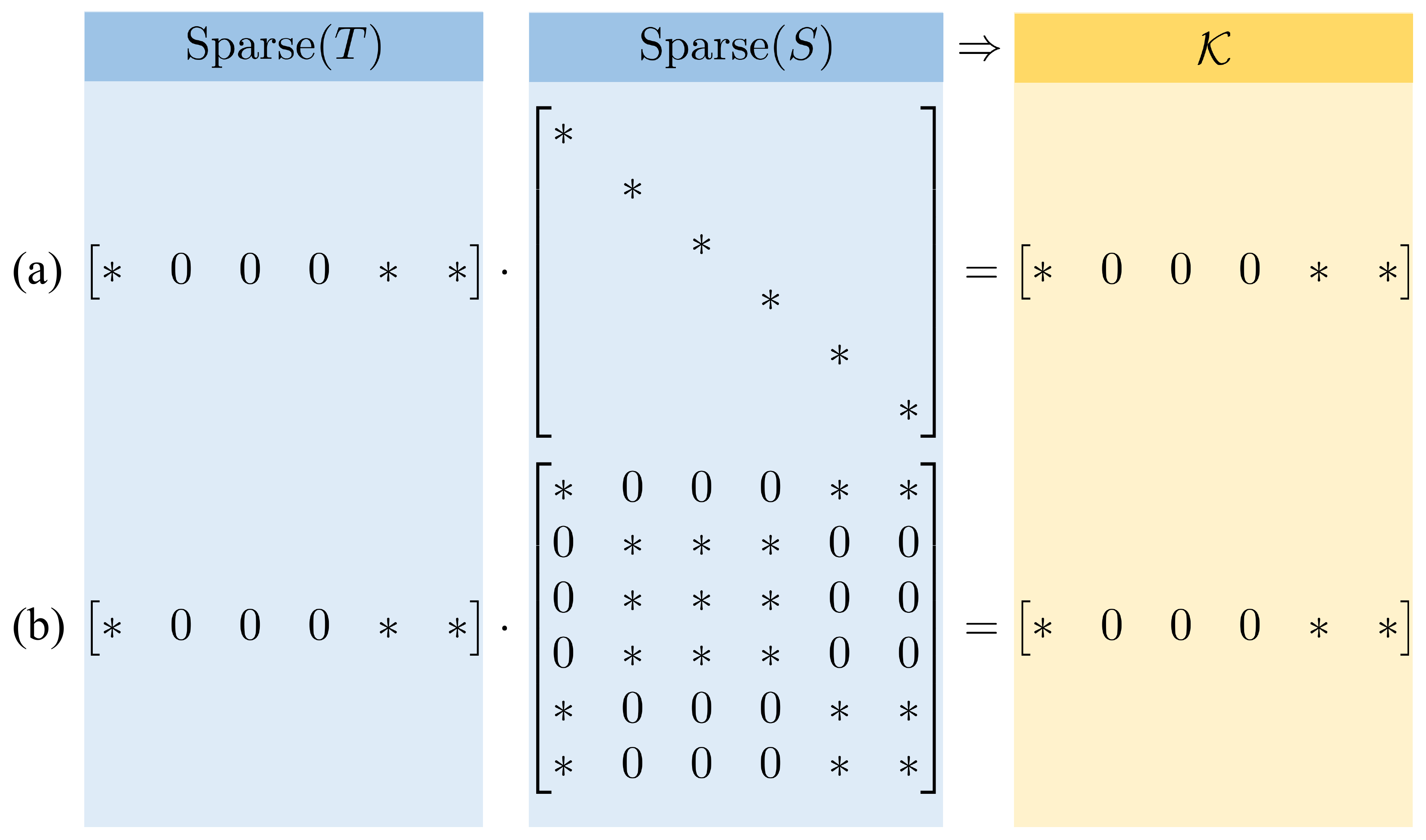}
	\caption{Two examples of possibles choices of $T$ and $S$ to guarantee sparsity invariance for a specific $\mathcal{K}$. (a) is to assume $\mathrm{Sparse}(S)$ to be diagonal, while (b) is a more general case.}
	\label{Fig:SparsityInvariance}
\end{figure}

Consider the equivalent formulation \eqref{Pr:P2}, where $ZX^{-1}\in \mathcal{K}$ is non-convex. The idea of sparsity invariance is to replace this non-convex constraint with separate constraints on $Z$ and $X$~\cite{furieri2019OnSep}. Precisely, we aim to characterize a binary matrix $T \in \{0,1\}^{1 \times 2n}$ and a symmetric binary matrix $S \in \{0,1\}^{2n\times 2n}$ such that
\begin{equation}\label{Eq:SparsityInvariance}
Z\in \mathrm{Sparse}(T)\; \mathrm{and}\; X\in\mathrm{Sparse}(S)  \Rightarrow ZX^{-1}\in \mathcal{K}.
\end{equation}
The above property is called \emph{sparsity invariance}; see Fig.~\ref{Fig:SparsityInvariance} for examples. If such $T$ and $S$ can be identified, then problem \eqref{Pr:P2} can be relaxed to the following convex optimization problem
\begin{equation} \label{Pr:P3}
\begin{aligned}
\min_{X,Y,Z} \quad &\mathrm{Tr}(QX)+\mathrm{Tr}(RY)\\
\mathrm{subject\;to}\quad & AX+XA^\tr-BZ-Z^\tr B^\tr+HH^\tr\preceq 0,\\
&\begin{bmatrix}
Y&Z\\Z^\tr&X
\end{bmatrix}\succeq 0,\;X\succ 0,\\
&Z\in \mathrm{Sparse}(T),\;X\in\mathrm{Sparse}(S).
\end{aligned}
\end{equation}
where any solution $(X,Y,Z)$ to \eqref{Pr:P3} is a suboptimal solution to \eqref{Pr:P2} due to the sparsity invariance property \eqref{Eq:SparsityInvariance}.

One key step in \eqref{Pr:P3} is to design two sparsity patterns $T$ and $S$ that satisfy \eqref{Eq:SparsityInvariance} for a given $\mathcal{K}$. A typical method, \eg, \cite{zheng2017convex}, is to assume $X$ to be block diagonal and $Z \in \mathcal{K}$; see Fig.~\ref{Fig:SparsityInvariance}a for illustration. Essentially, this assumption requires that the closed-loop system admits a block-diagonal Lyapunov function. However, this requirement might be too restrictive, and it fails to return feasible solutions in some instances. A full characterization of $T$ and $S$ was presented in \cite{furieri2019OnSep}. Following the strategy therein, we first choose a binary matrix $T$ to induce the same sparsity pattern $\mathcal{K}$, and then use the following two-step procedure to derive the optimal choice of $S$.

\begin{step}
	For every $i$, $j\in \{1,2,\ldots,n\}$, set
	\begin{equation}
	S_{ij}=
	\begin{cases}
	0,\;\mathrm{if}\,\exists k\in\{1,\ldots,m\}\, \mathrm{s.t.} \,T_{kj}=0, T_{ki}=1,\\
	1,\;\mathrm{otherwise}.
	\end{cases}
	\end{equation}
\end{step}

\begin{step}
	For every $i$, $j\in \{1,2,\ldots,n\}$, set
	\begin{equation}
	S_{ij}^*=
	\begin{cases}
		1,\;\mathrm{if}\,S_{ij}=S_{ji}=1,\\
		0,\;\mathrm{otherwise}.
	\end{cases}
	\end{equation}
\end{step}
It is shown in \cite{furieri2019OnSep} that $S^*$ is an optimal choice that maximizes the number of non-zero entries in $S$ for a given $T$, while satisfying the sparsity invariance property.

After obtaining the optimal choice of $S$ and $T$, problem \eqref{Pr:P3} is convex and ready to be solved by existing conic solvers (\eg, Mosek). One structured optimal controller is recovered as $K=ZX^{-1} \in \mathcal{K}$, which naturally satisfies the structural constraint. Note that the resulting optimal controller is in general a suboptimal solution to the original problem \eqref{Pr:P2}, since $Z$ and $X$ can be dense and still satisfy $ZX^{-1}\in \mathcal{K}$. Remarkably, however, our numerical experiments confirm that the suboptimal solution returned by solving \eqref{Pr:P3} often has satisfactory performance. In addition, if the entire traffic state is observable to the CAV, \ie, the structural constraint is removed, problem \eqref{Pr:P2} becomes a convex problem which returns a globally optimal controller.

\begin{figure}[t]
	\centering
	\includegraphics[scale=0.31]{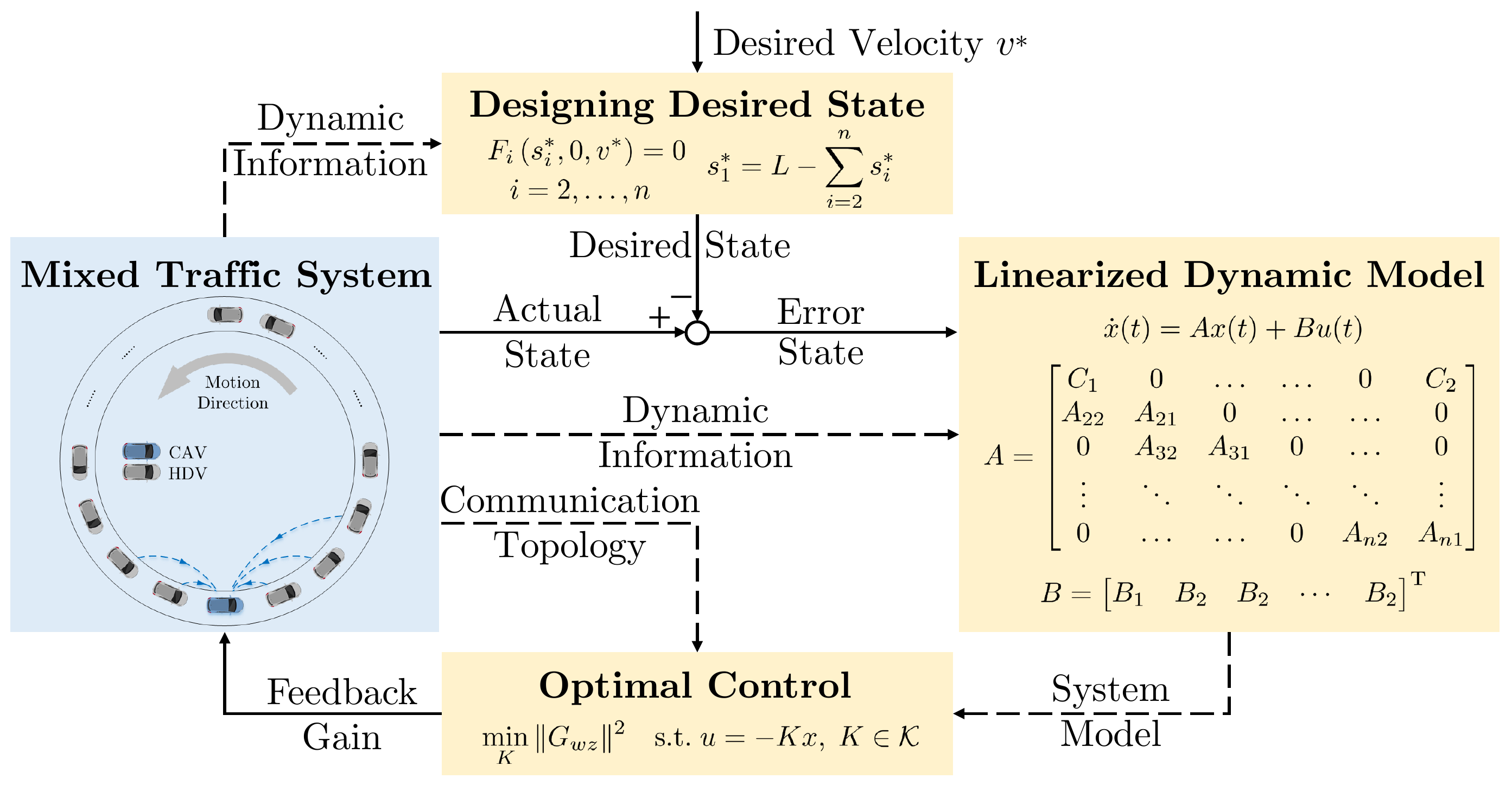}
	\caption{Process of designing an optimal controller for the mixed traffic system. Given a mixed traffic system, a pre-specified communication topology of the CAV and a desired equilibrium velocity, following this process allows the controller to calculate a static structured feedback gain $K$. The control input of the CAV at time $t$ is then computed by $u(t)=-Kx(t)$.}
	\label{Fig:DesignProcess}
\end{figure}

\subsection{Reachability Analysis: Desired Traffic Velocity}
\label{Sec:4c}

Solving problem \eqref{Pr:Original} offers a stabilizing feedback gain $K$ under structural constraints. Still, whether the desired traffic state can be reached remains unclear due to the existence of the zero eigenvalue. The mixed traffic system might converge to a different equilibrium traffic velocity $\hat{v}^{*}$, instead of the desired one $v^{*} $. We call this problem \emph{reachability analysis}.

The pre-specified velocity $v^{*} $ characterizes the desired state of the mixed traffic system, since each HDV has a corresponding equilibrium spacing  $ s_{i}^{*} $, determined by \eqref{Eq:Equilibrium}. Unlike HDVs, the equilibrium spacing of the CAV, \ie, $s_1^*$, can be designed separately. Consider the implementation of the controller, given by
\begin{equation}\label{Eq:ControlInput}
\begin{aligned}
u(t)=&-Kx(t)\\
=&-k_{11}\left(s_1(t)-s_1^*\right)-k_{12}\left(v_1(t)-v^*\right)\\
&-\sum\nolimits_{(i,1)\in\mathcal{E}^c,i\neq1 }\left(k_{i1}(s_i(t)-s_i^*)-k_{i2}(v_i(t)-v^*)\right).
\end{aligned}
\end{equation}
If $s_1^*$ is chosen carefully, the mixed traffic flow can be steered towards the exact velocity $v^*$. This is summarized in the following theorem.
\begin{theorem}[Reachability]
	\label{Theorem:Reachability}
	Consider the mixed traffic system with one CAV and $n-1$ heterogeneous HDVs given by \eqref{Eq:SystemModel}. Suppose that a static feedback gain is found by \eqref{Pr:Original} and the coefficient matrix in~\eqref{Eq:FinalLinearEquations} is non-singular. Then, the traffic flow can maintain stability at velocity $v^*$ if and only if the desired spacing of the CAV satisfies
	\begin{equation}\label{Eq:Rechability}
	s_1^*=L-\sum_{i=2}^{n}s_i^*,
	\end{equation}
	with $s_i^*$, $i=2,\ldots,n$ determined by \eqref{Eq:Equilibrium}.
\end{theorem}

\begin{IEEEproof}
	Since problem \eqref{Pr:Original} has a solution $K\in \mathcal{K}$, the mixed traffic system \eqref{Eq:SystemModel} is stable by controlling the CAV with input $u(t)=-Kx(t)$. Suppose that system \eqref{Eq:SystemModel} reaches its equilibrium state at $t_f$, then we have $\dot x (t_f)=0$, yielding $u(t_f )=0$, and $\dot{ \tilde{s}}_i (t_f )=0,\,\dot {\tilde{v}}_i (t_f )=0,\,i=1,\ldots,n$. According to \eqref{Eq:LinearHDVModel} and \eqref{Eq:LinearCAVModel}, we observe
	\begin{equation}\label{Eq:FinalHDV}
	\tilde{s}_i(t_f)=\frac{\alpha_{i2}-\alpha_{i3}}{\alpha_{i1}}v_e,\,\tilde{v}_i(t_f)=v_e,\,i=2,\ldots,n,
	\end{equation}
	and
	\begin{equation}\label{Eq:FinalCAV}
	\tilde{s}_1(t_f)=s_e,\,\tilde{v}_1(t_f)=v_e,
	\end{equation}
	where $s_e,v_e$ are two constants. Note that $\tilde{s}_i(t_f)=0,\,\tilde{v}_i(t_f)=0,\,i=1,\ldots,n$ is equivalent to the condition that the traffic system reaches the pre-specified velocity $v^*$. Due to the existence of the zero eigenvalue (as proved in Section \ref{Sec:3}),
	\begin{equation}\label{Eq:SpacingConstant}
	\sum_{i=1}^{n}s_i(t)=\sum_{i=1}^{n}\left(\tilde{s}_i(t)+s_i^*\right)=L,
	\end{equation}
	holds for all $t\geq 0$. Besides, $u(t_f)=0$ yields
	\begin{equation}\label{Eq:FinalInput}
	Kx(t_f)=0.
	\end{equation}
	Substituting \eqref{Eq:FinalHDV} and \eqref{Eq:FinalCAV} into \eqref{Eq:SpacingConstant} and \eqref{Eq:FinalInput}, we have
	\begin{equation}\label{Eq:FinalLinearEquations}
	\begin{cases}
	s_e+\left(\sum_{i=2}^{n} \frac{\alpha_{i2}-\alpha_{i3}}{\alpha_{i1}}\right)v_e=L-s_1^*-\sum_{i=2}^{n}s_i^*,\\
	k_{11}s_e+\left(\sum_{(i,1)\in\mathcal{E}^c,i\neq1 }\left(k_{i1} \frac{\alpha_{i2}-\alpha_{i3}}{\alpha_{i1}} +k_{i2}\right)+k_{12} \right)v_e=0.
	\end{cases}
	\end{equation}
	
	To show sufficiency, we assume that \eqref{Eq:Rechability} holds. The linear equation system \eqref{Eq:FinalLinearEquations} then becomes a homogeneous one, which has a unique solution most generally: $s_e=0,\,v_e=0$. This result indicates $\tilde{s}_i(t_f)=0,\,\tilde{v}_i(t_f)=0,\,i=1,\ldots,n$. Thus, the traffic system can reach the desired traffic velocity $v^*$.
	
	To show necessity, we assume $\tilde{s}_i(t_f)=0,\,\tilde{v}_i(t_f)=0,\,i=1,\ldots,n$. Applying \eqref{Eq:FinalHDV} and \eqref{Eq:FinalCAV}, we see $s_e=0,\,v_e=0$. Substituting it into \eqref{Eq:FinalLinearEquations} yields $L-s_1^*-\sum_{i=2}^{n}s_i^* =0$, which proves that \eqref{Eq:Rechability} holds.
	
\end{IEEEproof}

\begin{remark}
	The whole process of obtaining a system-level optimal controller is illustrated in Fig.~\ref{Fig:DesignProcess}, including the problem formulation, the numerical solution, and the design of the desired system state. This process offers a static linear feedback gain $K$, which allows the CAV to utilize local available information and achieve an optimal performance for the entire traffic flow. Note that our strategy requires the explicit model of the HDVs' dynamics, which is also assumed in other recent works; see, \emph{e.g.}, \cite{jin2017optimal,wu2018Stabilizing,cui2017stabilizing,vaio2019Cooperative}. In practice, our proposed controller can be combined with existing real-time algorithms, \emph{e.g.}, \cite{pandita2013preceding,jin2018connected}, for estimating HDVs' dynamics based on historical data of vehicle trajectories.
	Nevertheless, considering possible model mismatch in HDVs' dynamics, it is of great significance to incorporate robustness analysis in the optimal controller synthesis in the future work.
\end{remark}

In design of ACC or CACC \cite{li2017dynamical}, the notion of the desired state also exists, which can be chosen arbitrarily since all the involved vehicles have autonomous capabilities. In contrast, in mixed traffic systems, one distinctive fact is that the HDVs cannot be directly controlled. Despite this fact, they can still be influenced indirectly through interactions among neighboring vehicles towards an equilibrium state $(s_i^*,v^*)$. In addition, although the desired state of the CAV, \ie, $s_1^*$, can be chosen separately, it should satisfy \eqref{Eq:Rechability} in order to stabilize the traffic flow at the desired velocity $v^*$. If $s_1^*\neq L-\sum_{i=2}^{n}s_i^*$, then it is easy to see from \eqref{Eq:FinalLinearEquations} that $s_e\neq 0,\,v_e\neq 0$, leading to $v_i (t_f )=v_f\neq v^*$, $i=1,\ldots,n$. This means that the traffic system is stabilized at another equilibrium velocity $v_f$ instead of the pre-specified one $v^*$.

Moreover, we observe that the pre-specified traffic velocity $v^*$ should also satisfy a certain constraint to make it reachable. This is summarized as follows.
\begin{corollary} [Maximum reachable traffic velocity]
	\label{Corollary:MaximumVelocity}
	Considering the mixed traffic system given by \eqref{Eq:SystemModel}, there exists a reachable range for the traffic velocity:
	\begin{equation}\label{Eq:VelocityRange}
	0\le v^*<v^*_{\max},
	\end{equation}
	where $v^*_{\max}$ denotes the maximum reachable traffic velocity.
\end{corollary}
\begin{IEEEproof}
	Since the spacing of the CAV must be positive, we have $s_1^*>0$. Substituting \eqref{Eq:Rechability} into $s_1^*>0$ yields
	\begin{equation}\label{Eq:PositiveSpacing}
	L-\sum_{i=2}^{n}s_i^*>0.
	\end{equation}
 Recall that the HDV equilibrium equation \eqref{Eq:Equilibrium} is an implicit function associating $v^*$ and $s_i^*$. 
In general, $s_i^*$ increases as $v^*$ grows up \cite{jin2017optimal}, as has been shown in Fig.~\ref{Fig:HDVRelationship}. Thus, we suppose there exist non-decreasing functions $E_i(\cdot)$, $i=2,\ldots,n$, such that $s_i^*=E_i (v^*)$. Then, \eqref{Eq:PositiveSpacing} can be converted to
	$
	\sum_{i=2}^{n} E_i(v^*)<L,
	$
	leading to a maximum reachable velocity $v_{\max}^*$, given by
	\begin{equation}\label{Eq:MaximumVelocity}
	\sum_{i=2}^{n} E_i(v_{\max}^*)=L.
	\end{equation}
	Accordingly, \eqref{Eq:VelocityRange} holds.
\end{IEEEproof}

\begin{remark}
	Although the stabilizability of the mixed traffic system is revealed in Section~\ref{Sec:3}, the equilibrium velocity $v^*$ still needs to be designed carefully for practical use. Corollary \ref{Corollary:MaximumVelocity} indicates that it is not possible to stabilize the mixed traffic flow at an arbitrary velocity $v^*$; instead, there exists a maximum reachable traffic velocity $v_{\max}^*$. According to \eqref{Eq:MaximumVelocity}, this upper bound has a certain relationship with the number of vehicles, HDV dynamics, and the ring-road circumference. In the homogeneous case, \eqref{Eq:MaximumVelocity} reduces to ($E_i=E,\,i=2,\ldots,n$)
	\begin{equation}
	v_{\max}^*=E^{-1}\left(\frac{L}{n-1}\right),
	\end{equation}
	which is consistent with that in \cite{zheng2018smoothing}. Considering the monotonicity of function $E$, we conclude that the maximum reachable traffic velocity decreases as the vehicle density grows up. In addition, note that a larger value of $v^*$ leads to a smaller value of the CAV's desired spacing $s_1^*$, which will increase the risk of rear-end collisions. In contrast, a smaller value of $v^*$ leaves a larger spacing between the CAV and the preceding vehicle, which may cause other vehicles to cut in. Therefore, we need to choose a moderate value of $v^*$ among the range of $0 \leq v^* < v_{\max}^*$ for practical use, based on the trade-off between vehicle safety, practical feasibility, and traffic efficiency.
\end{remark}

\section{Numerical Experiments}
\label{Sec:5}
Our main results are obtained using a linearized model of the mixed traffic system. We evaluate their effectiveness in the presence of nonlinearities arising in the car-following dynamics \eqref{Eq:HDVModel}. In this section, we conduct three types of simulation experiments to validate our theoretical results based on a realistic nonlinear HDV model. All the experiments are carried out in MATLAB.

\subsection{Experimental Setup}

We consider a ring road with circumference $L=400~m$ containing 19 HDVs and one CAV. Vehicle no.1 is the CAV and the penetration rate is only $5\%$ in this setup. The information from five vehicles ahead and five vehicles behind is available to the CAV for feedback control. For the parameters in the performance output \eqref{Eq:Output}, we choose $\gamma_s=0.03$, $\gamma_v=0.15$, $\gamma_u=1$. Based on the approach in Section \ref{Sec:4}, a structured linear feedback gain $K$ is obtained using Mosek.

In our simulations, a nonlinear optimal velocity model (OVM) \cite{bando1995dynamical} is used to describe the car-following dynamics of HDVs. The specific expression for \eqref{Eq:HDVModel} becomes
\begin{equation} \label{Eq:OVM}
F_i(\cdot)=\alpha_i\left(V_i\left(s_i(t)\right)-v_i(t)\right)+\beta_i\dot{s}_i(t),
\end{equation}
where $\alpha_i$, $\beta_i$ are sensitivity coefficients, and $V_i (s)$ is the spacing-dependent desired velocity of driver $i$, typically given by a continuous piecewise function
\begin{equation}
\label{Eq:OVMSpacingPolicy}
V_i(s)=\begin{cases}
0, &s\le s_{i,\mathrm{st}},\\
f_{i,v}(s), &s_{i,\mathrm{st}}<s<s_{i,\mathrm{go}},\\
v_{i,\max}, &s\ge s_{i,\mathrm{go}},
\end{cases}
\end{equation}
In \eqref{Eq:OVMSpacingPolicy}, we choose a nonlinear form of $f_{i,v}(s)$ \cite{jin2017optimal}
\begin{equation}
f_{i,v}(s)=\frac{v_{i,\max }}{2}\left(1-\cos (\pi
\frac{s-s_{i,\mathrm{st}}}{s_{i,\mathrm{go}}-s_{i,\mathrm{st}}})\right).
\end{equation}

Due to the heterogeneity setup and motivated by \cite{jin2017optimal}, we set the parameters as follows: $\alpha_i=0.6+U[-0.1,0.1]$, $\beta_i=0.9+U[-0.1,0.1]$, $s_{i,\mathrm{go}}=35+U[-5,5]$, where $U[\cdot]$ denotes the uniform distribution. The rest of parameters are set as $v_{i,\max}=30$, $s_{i,\mathrm{st}}=5$. As stated in \cite{zheng2018smoothing}, this set of values captures a stop-and-go car-following behavior of HDVs in a ring road. It is not difficult to verify that the stabilizability condition \eqref{Eq:ControllabilityTheorem} is satisfied. Accordingly, Theorem \ref{Theorem:Stabilizability} guarantees that the mixed traffic system is stabilizable by a single CAV around an equilibrium traffic state.

Note that our experimental setup is similar to that in \cite{stern2018dissipation,sugiyama2008traffic}, but has a larger inter-vehicle distance, which leads to higher traffic velocity. Compared to \cite{stern2018dissipation,sugiyama2008traffic}, our setup captures more common scenarios on highways, and meanwhile reproduces the experimental observations in \cite{stern2018dissipation,sugiyama2008traffic}; see Sections \ref{Sec:StabilizingTraffic} and \ref{Sec:DissipateWaves}. To avoid crashes, we also assume that all the vehicles are equipped with a standard automatic emergency braking system, described as follows
\begin{equation}
	\dot{v}(t)=a_{\min },\,\mathrm{if}\,\frac{v_{i}^{2}(t)-v_{i-1}^{2}(t)}{2s_{i}(t)}\ge \vert a_{\min }\vert,
\end{equation}
where the maximum acceleration and deceleration rate of each vehicle are set to $a_{\max}=2 \,\mathrm{m/s^2}$, $a_{\min}=-5 \,\mathrm{m/s^2}$, respectively.

\subsection{Stabilizing Mixed Traffic Flow}
\label{Sec:StabilizingTraffic}

\begin{figure}[t]
	\centering
	\subfigure[]
	{\label{Fig:ExperimentAa}
		\includegraphics[scale=0.35]{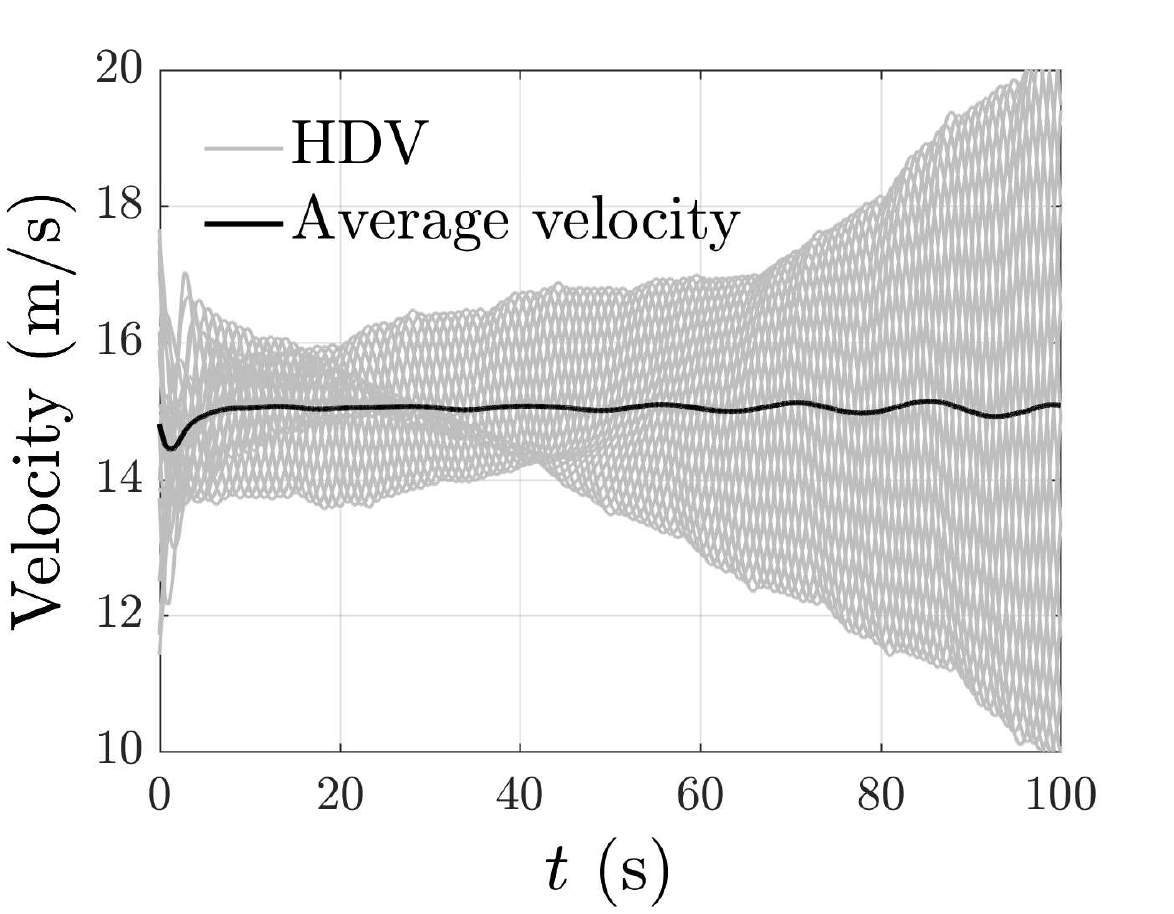}}
	\subfigure[]
	{\label{Fig:ExperimentAb}
		\includegraphics[scale=0.35]{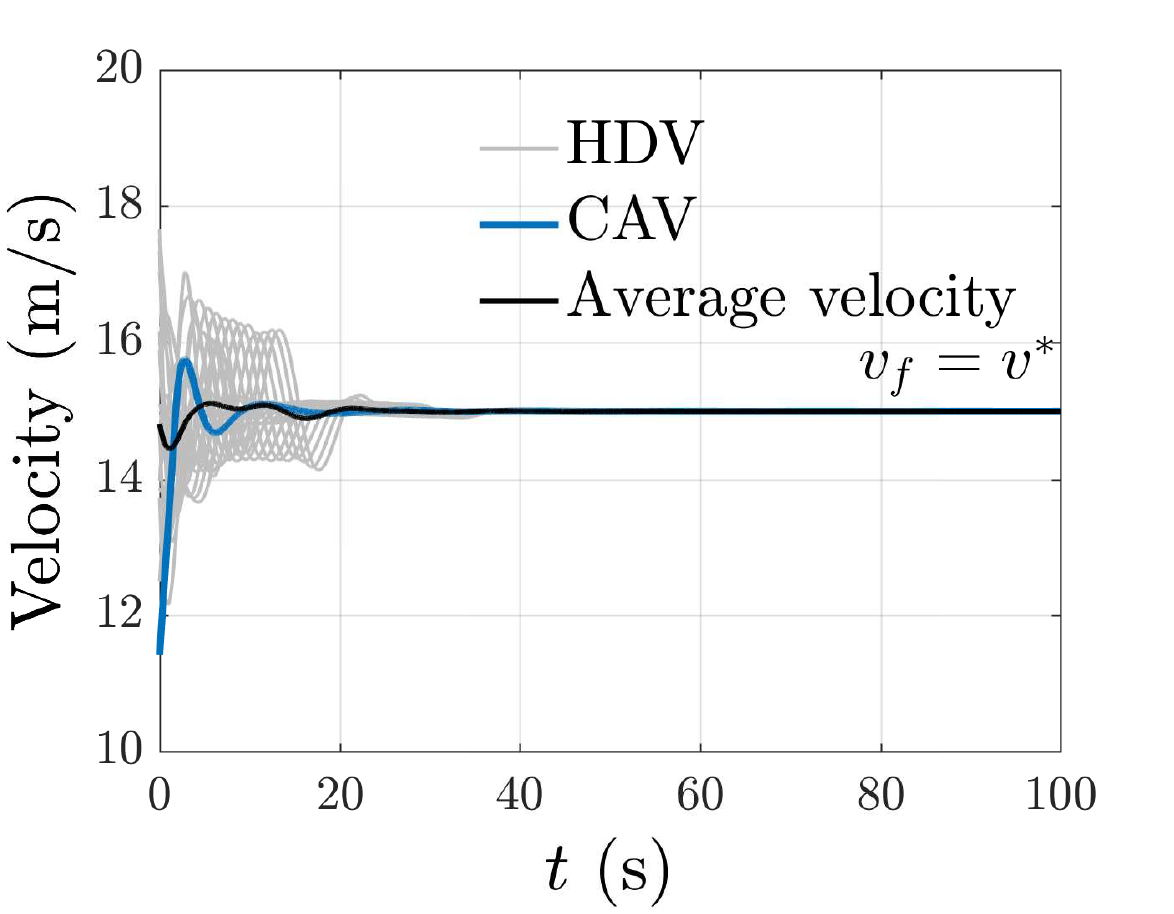}}
	\subfigure[]
	{\label{Fig:ExperimentAc}
		\includegraphics[scale=0.35]{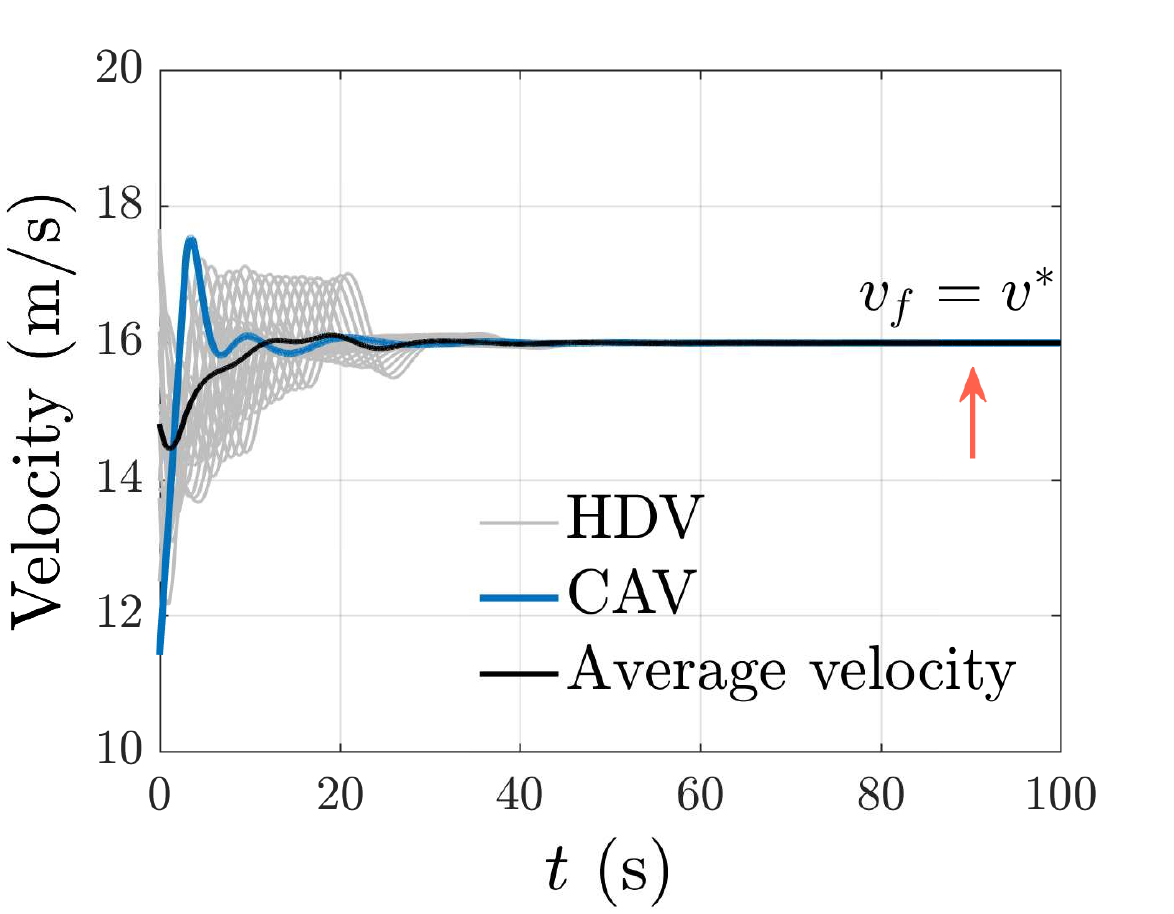}}
	\subfigure[]
	{\label{Fig:ExperimentAd}
		\includegraphics[scale=0.35]{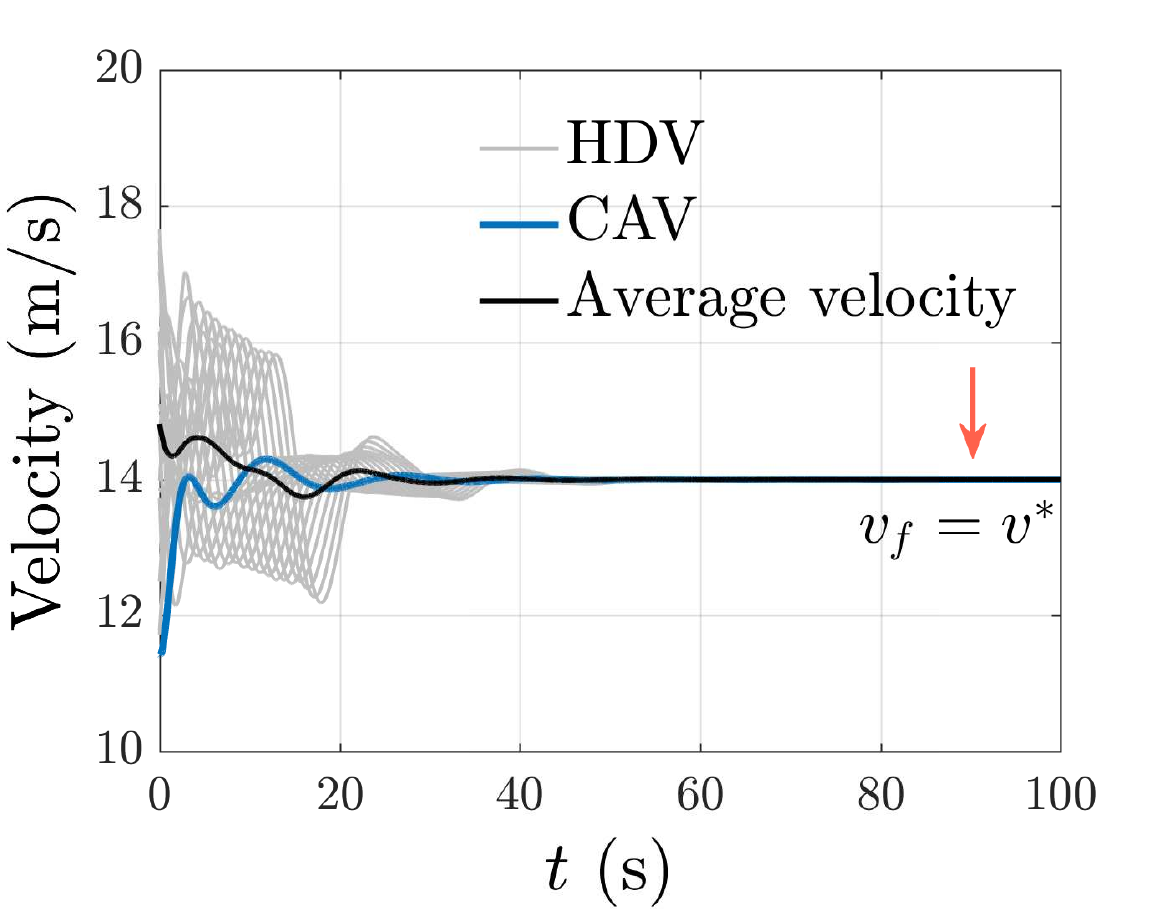}}
	\subfigure[]
	{\label{Fig:ExperimentAe}
		\includegraphics[scale=0.35]{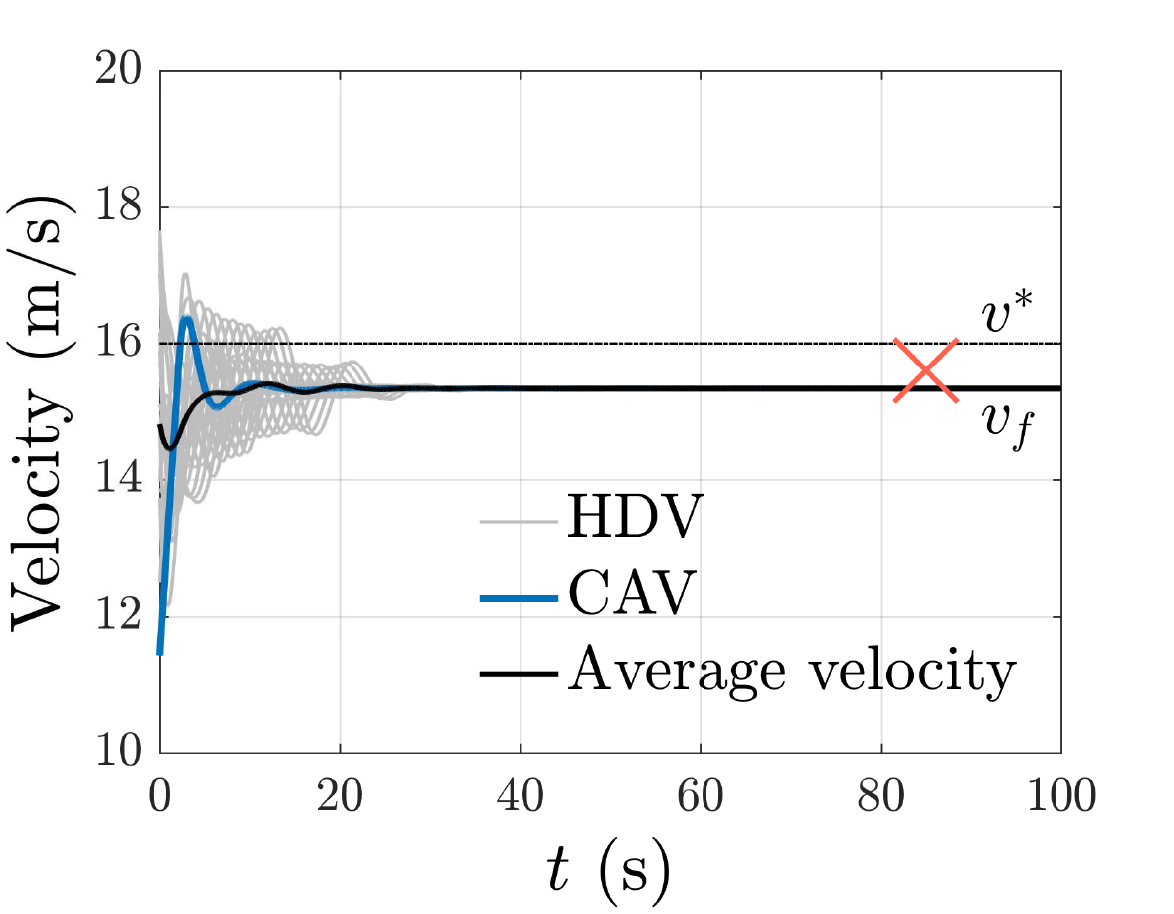}}
	\subfigure[]
	{\label{Fig:ExperimentAf}
		\includegraphics[scale=0.35]{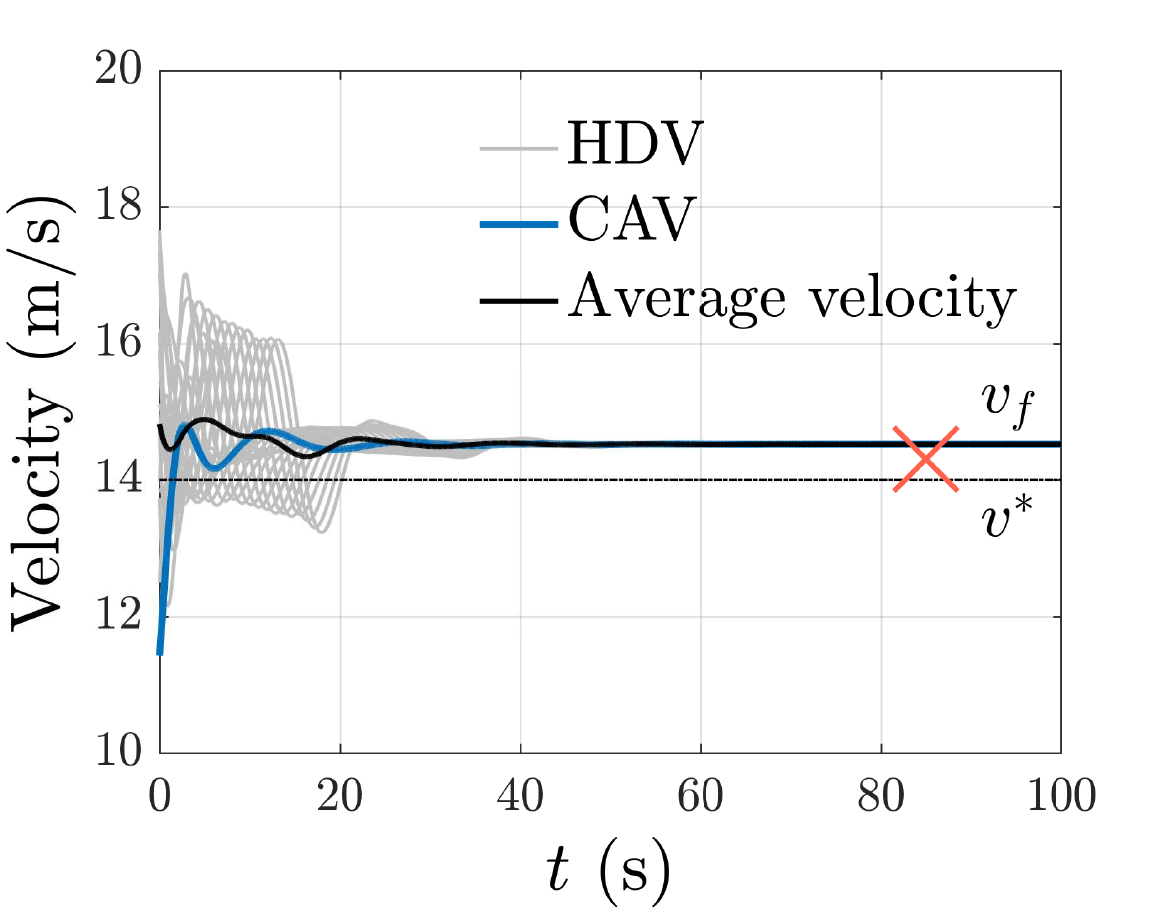}}
	\caption{Velocity profile of each vehicle (Experiment A). (a) All the vehicles are HDVs. (b)---(f) One vehicle is CAV with the proposed method. The desired equilibrium velocity $v^*$ is $15\,\mathrm{m/s}$, $16\,\mathrm{m/s}$, $14\,\mathrm{m/s}$ in (b), (c)(e), (d)(f), respectively. In (b)(c)(d) the value of $s_1^*$ is determined according to \eqref{Eq:Rechability}, whereas in (e)(f) \eqref{Eq:Rechability} is not satisfied.}
	\label{Fig:ExperimentA}
\end{figure}

Experiment A aims to examine the analytical results with respect to stabilizability and reachability, as shown in Theorems~\ref{Theorem:Stabilizability} and \ref{Theorem:Reachability} respectively. Firstly, we show that a typical nonlinear heterogeneous traffic flow can be stabilized by a single CAV. At the beginning of the simulation, all the vehicles are distributed randomly on the road, with the initial velocity following the distribution $15+U[-4,4] \,\mathrm{m/s}$. When all the vehicles are under human control, it is clearly observed that multiple perturbations arise inside the traffic flow, and they are amplified gradually, inducing a traffic wave propagating against the traffic flow (Fig.~\ref{Fig:ExperimentAa}). In contrast, if there is one CAV using the proposed control method, the traffic flow can be stabilized to the original average velocity $15\,\mathrm{m/s}$ within a short time (Fig.~\ref{Fig:ExperimentAb}). Moreover, by adjusting the desired equilibrium velocity $v^*$ (within the range as stated in Corollary \ref{Corollary:MaximumVelocity}) and the corresponding desired spacing $s_1^*$ according to Theorem~\ref{Theorem:Reachability}, the CAV shows its ability to steer the entire traffic flow towards a higher or lower velocity, through influencing other vehicles (Fig.~\ref{Fig:ExperimentAc}-\ref{Fig:ExperimentAd}).

Next we show that the desired state of the CAV needs to be designed carefully, which validates the reachability analysis. Following the same simulation setup as that in Fig.~\ref{Fig:ExperimentAc}-\ref{Fig:ExperimentAd}, we change the value of the desired spacing of the CAV, \ie, $s_1^*$, without satisfying \eqref{Eq:Rechability} in Theorem \ref{Theorem:Stabilizability}. Note that the feedback gain $K$ remains the same. As illustrated in Fig.~\ref{Fig:ExperimentAe}-\ref{Fig:ExperimentAf}, although the mixed traffic system is still stabilized by the CAV, the system final velocity $v_f$ is not equal to the desired velocity $v^*$. Due to the nonlinearity of traffic flow, the deviation becomes unpredictable and cannot be derived from \eqref{Eq:FinalLinearEquations}. Conversely, if $s_1^*$ is designed based on Theorem \ref{Theorem:Reachability}, Fig.~\ref{Fig:ExperimentAc}-\ref{Fig:ExperimentAd} exhibit that the traffic flow can be controlled precisely to the pre-specified traffic velocity. This result confirms the statement in Theorem~\ref{Theorem:Reachability}.

\subsection{Dissipating Stop-and-Go Waves}
\label{Sec:DissipateWaves}

\begin{figure}[t]
	\centering
	\subfigure[]
	{\hspace{-4mm}
		\label{Fig:ExperimentBa}
		\includegraphics[scale=0.34]{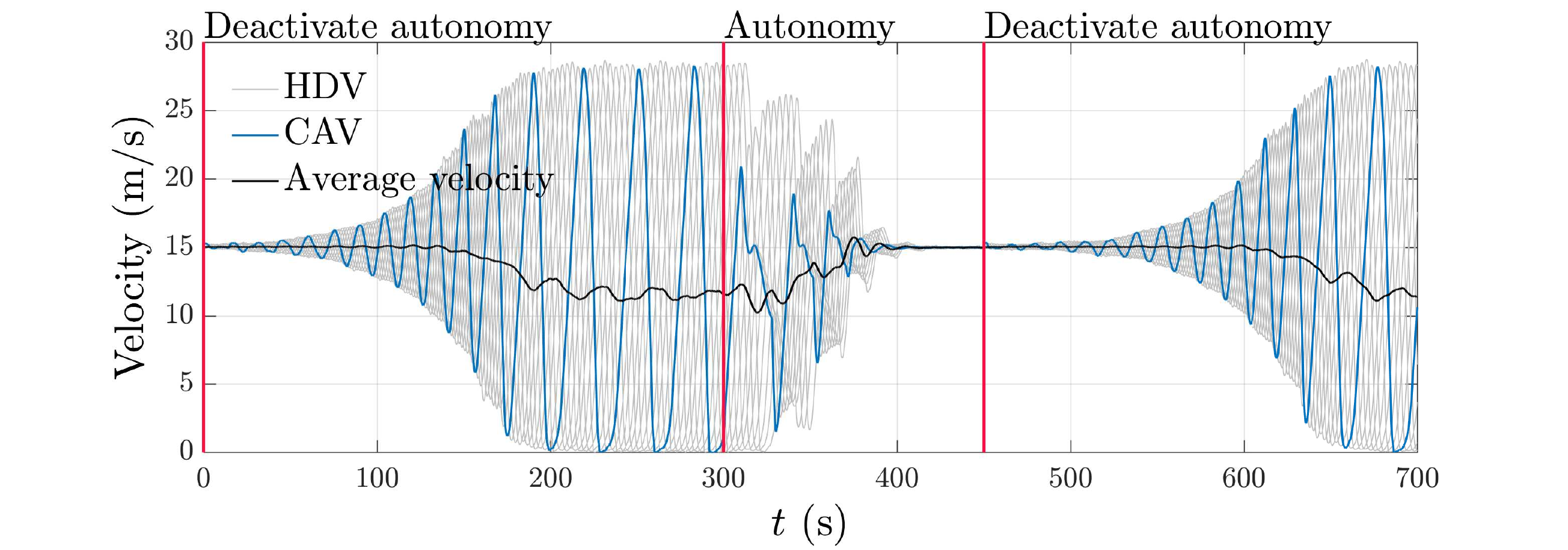}}
	\subfigure[]
	{\hspace{-4mm}
		\label{Fig:ExperimentBb}
		\includegraphics[scale=0.34]{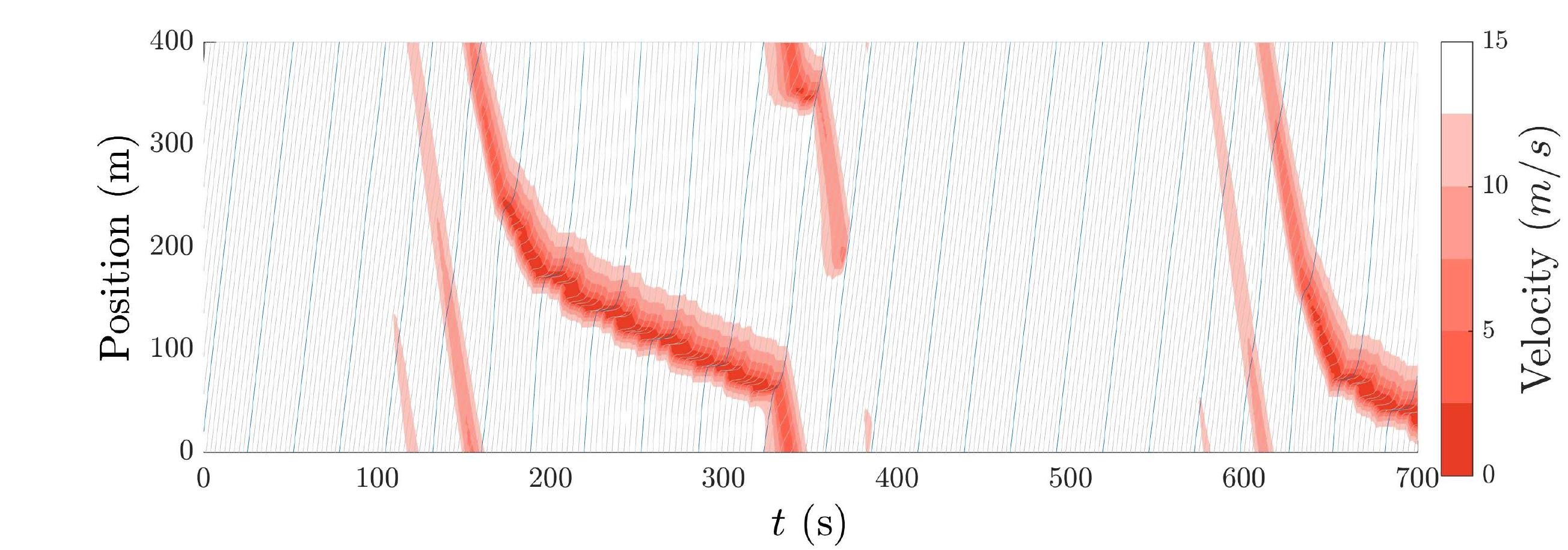}}
	\caption{Velocity profile and trajectory of the 1st, 3rd, 5th, $\ldots$, 19th vehicle (Experiment B). In (b), the darker the color, the lower the velocity. From $t=0\,\mathrm{s}$ to $t=300\,\mathrm{s}$ and from $t=450\,\mathrm{s}$ to $t=700\,\mathrm{s}$, the proposed controller does not work and all the vehicles are human-driven, while from $t=300\,\mathrm{s}$ to $t=450\,\mathrm{s}$, the proposed controller at vehicle no.1 is activated.}
	\label{Fig:ExperimentB}
\end{figure}

Experiment B is conducted to test the controller's ability to dissipate stop-and-go waves. A random noise following the normal distribution, $N(0,0.2)$, is added to the acceleration signal of each vehicle. This corresponds to the traffic situations where small perturbations are generated naturally inside the traffic flow. The perturbations may lead to traffic congestion, corresponding to traffic jams without bottleneck \cite{sugiyama2008traffic}. Our setup of random noise is consistent with that in the learning-based work \cite{wu2017flow}, and aims to reproduce the empirical observations in real-world experiments \cite{stern2018dissipation,sugiyama2008traffic} through our numerical studies.

The results are demonstrated in Figure \ref{Fig:ExperimentB}. At the beginning, the controller is deactivated, which means that vehicle no.1 behaves as an HDV, \ie, utilizes the OVM model \eqref{Eq:OVM} to determine its car-following behavior. This corresponds to the scenario where all the vehicles are controlled by human drivers. From $t=0\,\mathrm{s}$ to $t=300\,\mathrm{s}$, a traffic wave is observed, which grows up gradually and finally leads to a stop-and-go wave. At $t=300\,\mathrm{s}$, the proposed controller is activated, which dissipates the stop-and-go wave and moves the traffic flow back to equilibrium during the following $100$ seconds. At $t=450\,\mathrm{s}$, the controller is deactivated again, and the wave reappears. This result implies that the proposed method enables the CAV to attenuate the inner perturbations of the traffic flow persistently, which may be induced by system noise. Without the CAV using the proposed controller, the traffic flow may easily get into the congestion pattern.

\subsection{Comparison with Existing Heuristic Strategies}

\begin{figure}[t]
	\centering
	\subfigure[]
	{\label{Fig:ExperimentC_1a}
		\includegraphics[scale=0.35]{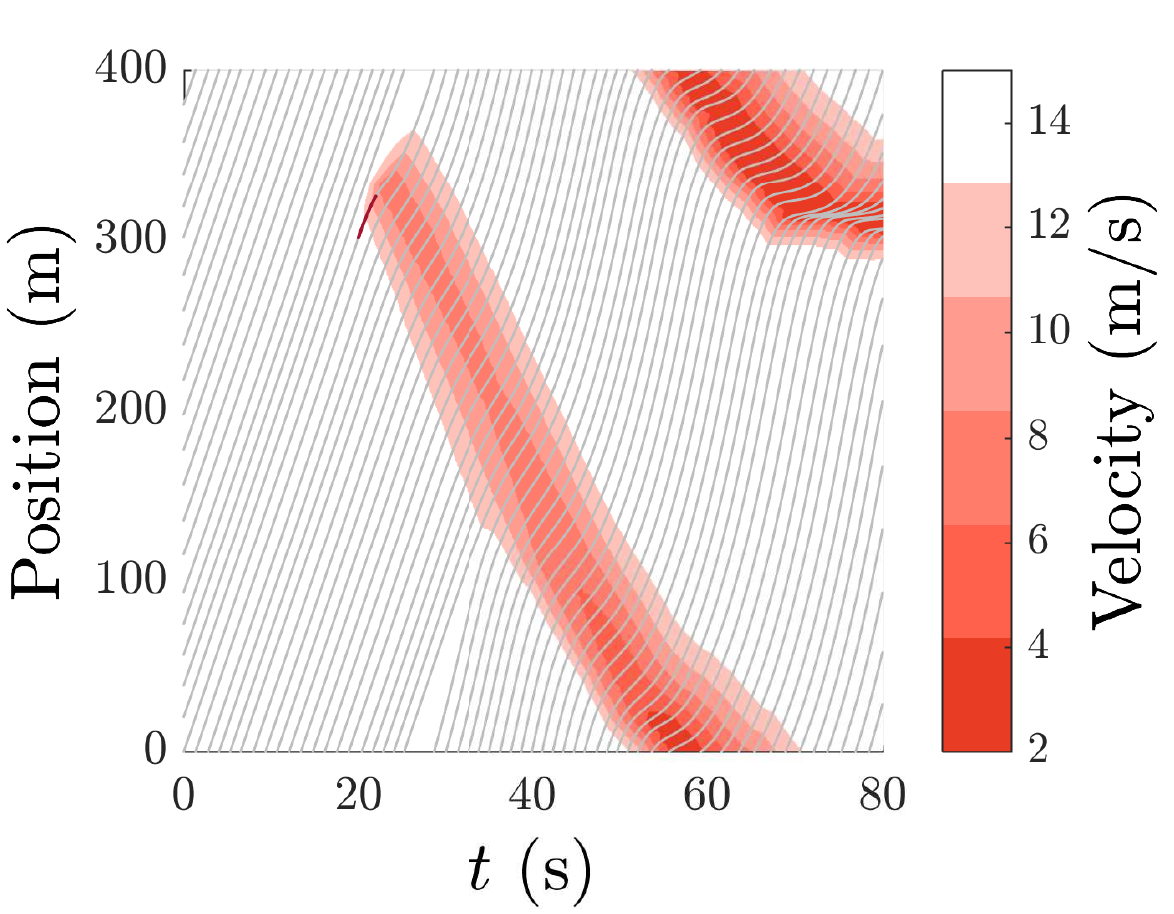}}
	\subfigure[]
	{\label{Fig:ExperimentC_1b}
		\includegraphics[scale=0.35]{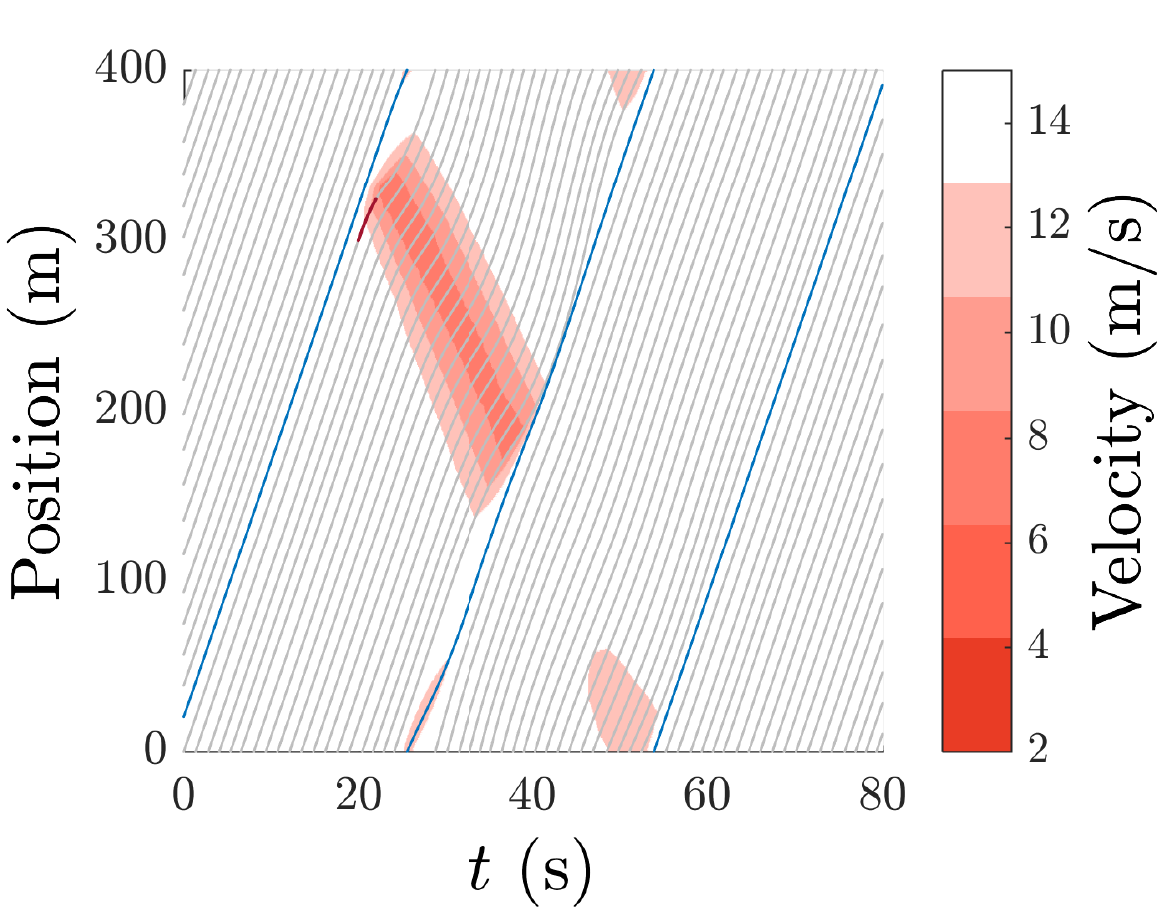}}
	\subfigure[]
	{\label{Fig:ExperimentC_1c}
		\includegraphics[scale=0.35]{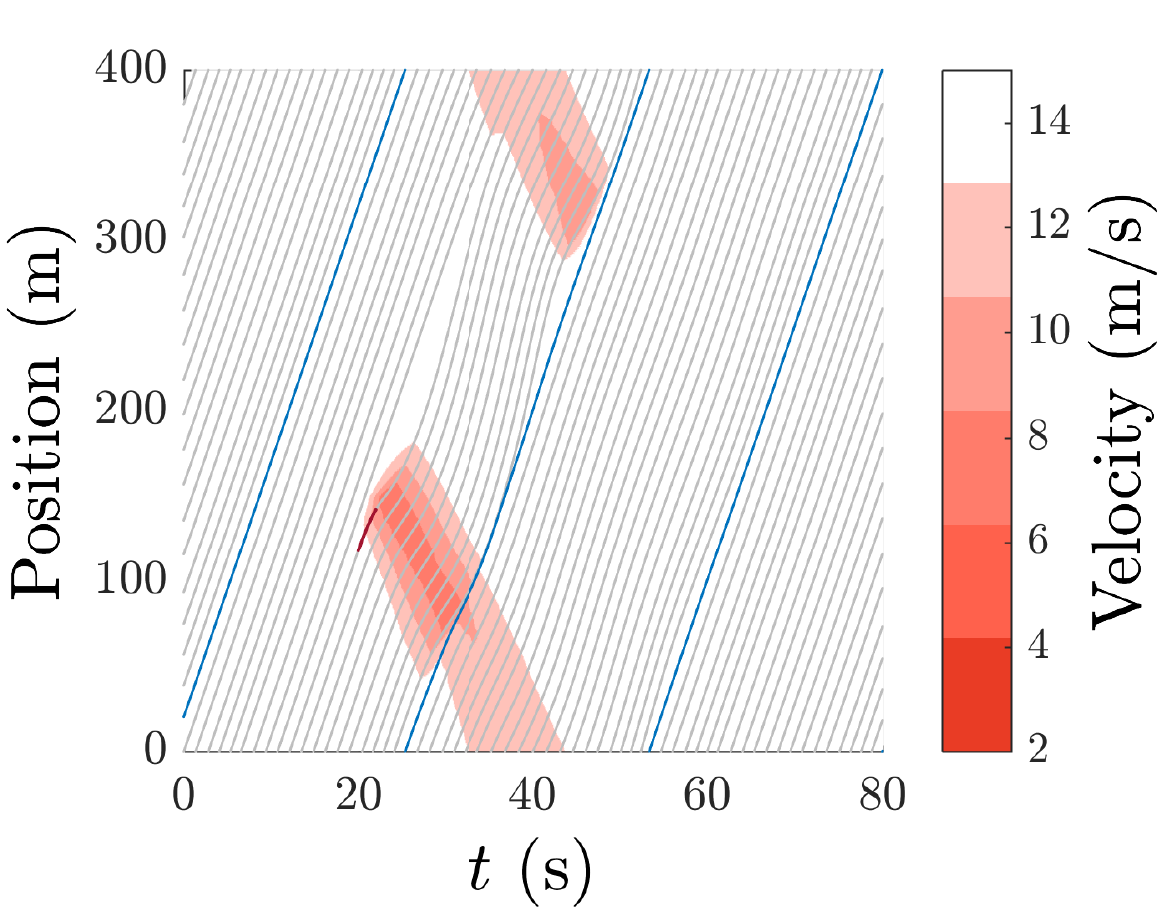}}
	\subfigure[]
	{\label{Fig:ExperimentC_1d}
		\includegraphics[scale=0.35]{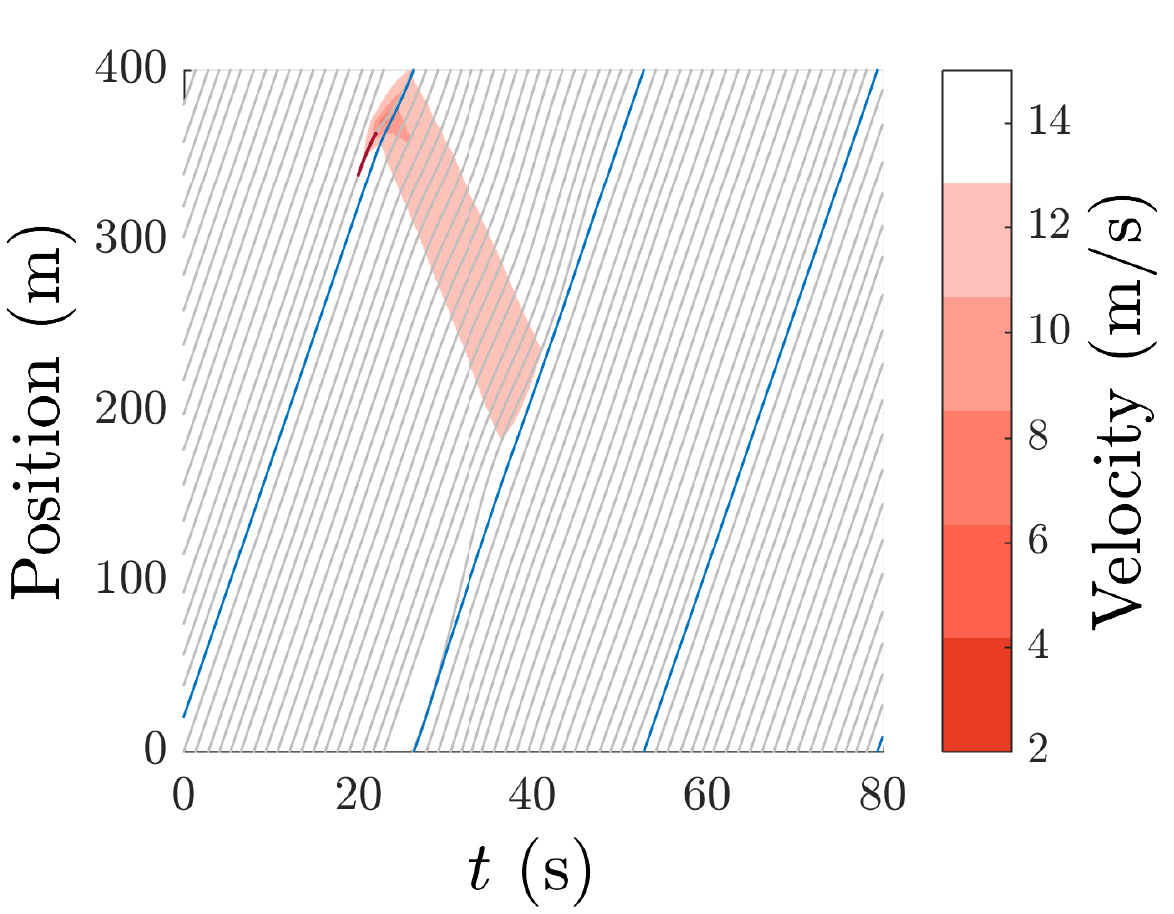}}
	\caption{Vehicle trajectories (Experiment C). (a) All the vehicles are human-driven. (b)-(d) correspond to the cases where vehicle 2,11,20 is under the perturbation, respectively.}
	\label{Fig:ExperimentC_1}
\end{figure}

\begin{figure}[t]
	\centering
	\subfigure[]
	{\label{Fig:ExperimentC_2a}
		\includegraphics[scale=0.35]{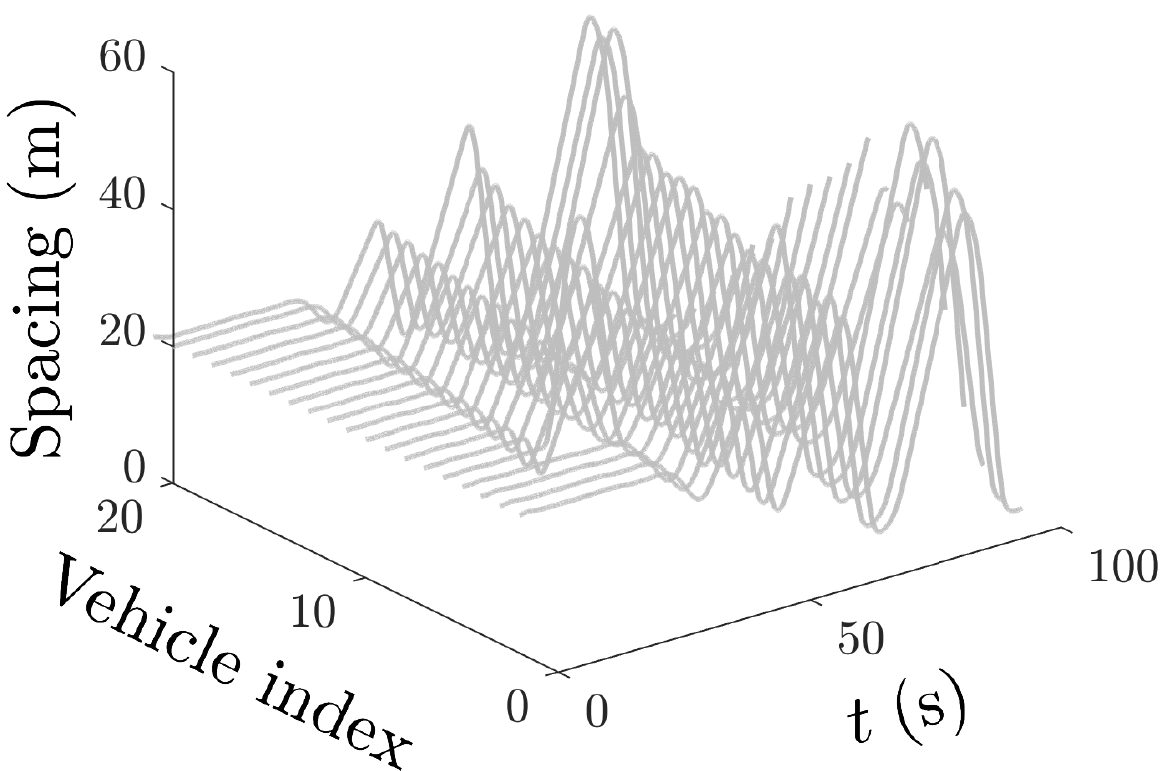}}
	\subfigure[]
	{\label{Fig:ExperimentC_2b}
		\includegraphics[scale=0.35]{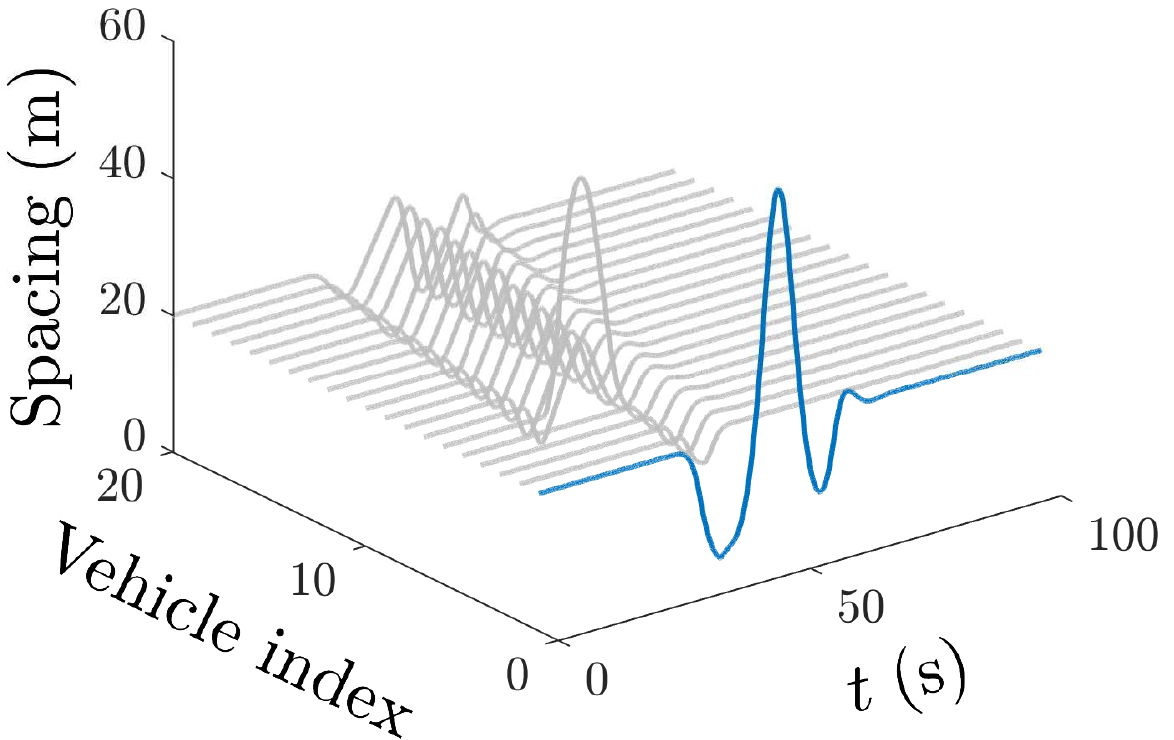}}
	\subfigure[]
	{\label{Fig:ExperimentC_2c}
		\includegraphics[scale=0.35]{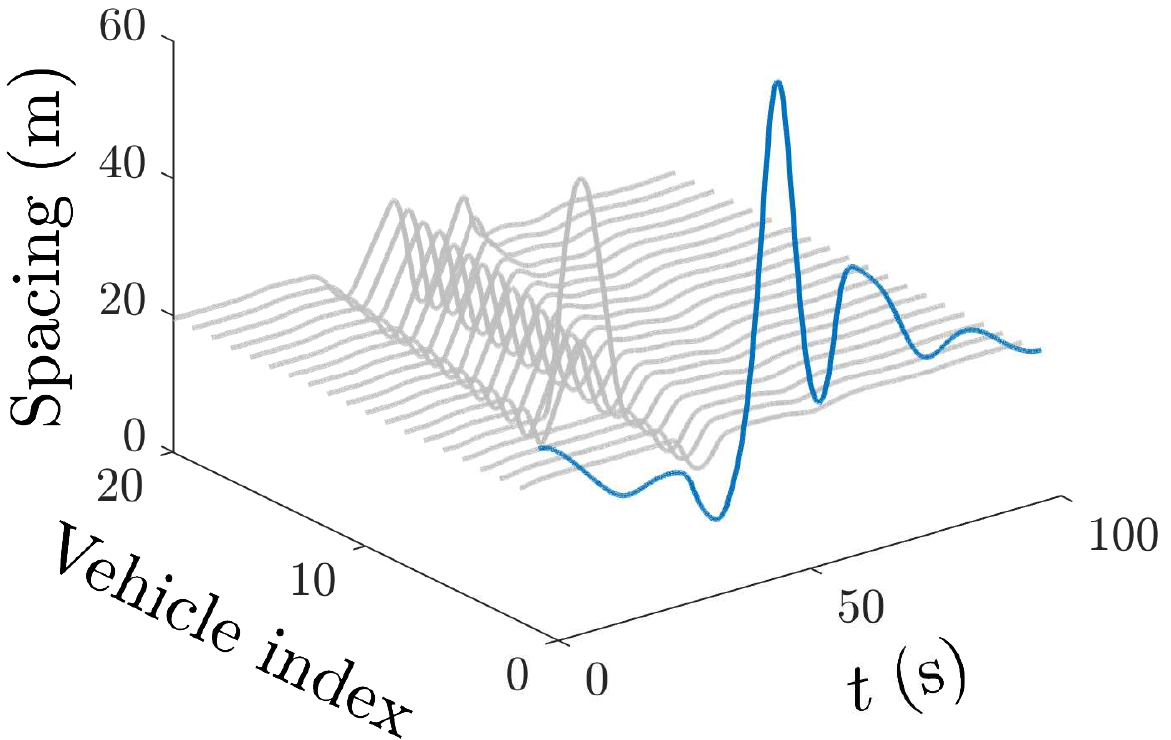}}
	\subfigure[]
	{\label{Fig:ExperimentC_2d}
		\includegraphics[scale=0.35]{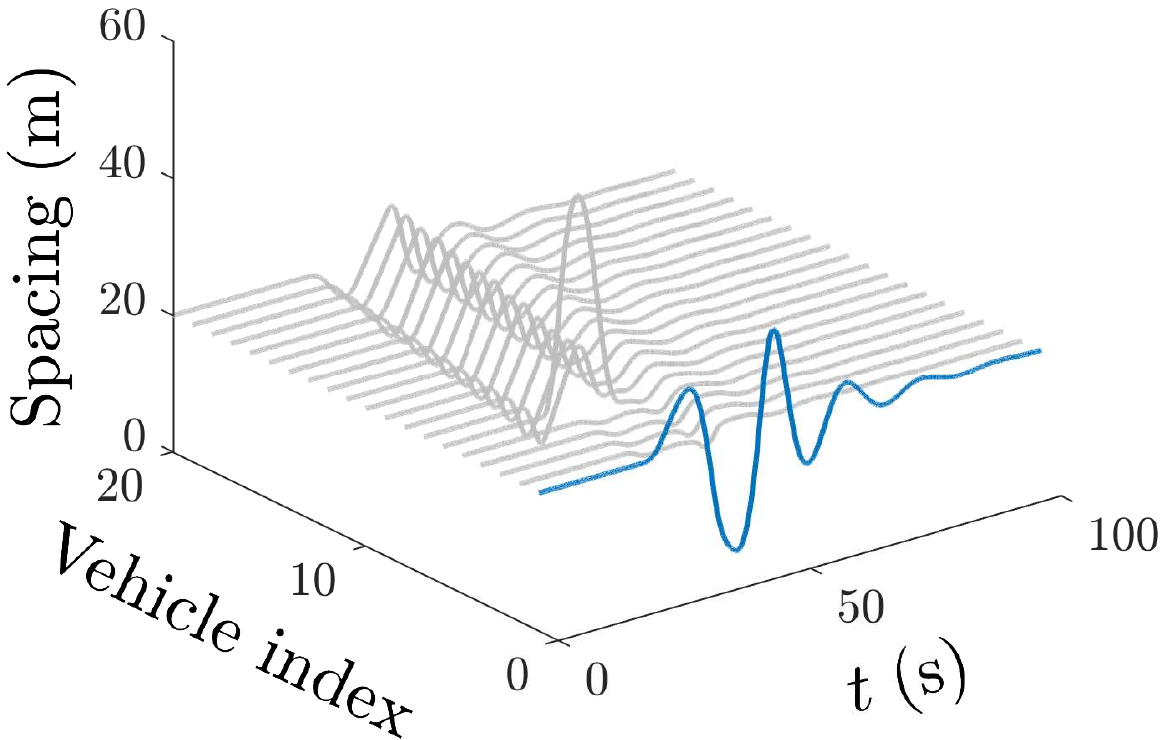}}
	\caption{Vehicle Spacing (Experiment C). Vehicle no.6 is under the perturbation. The grey line and the blue line denote the spacing of HDV and CAV, respectively. 
		(a) All the vehicles are HDVs. (b)-(d) correspond to the cases where the CAV is using FollowerStopper, PI with Saturation and our optimal control strategy, respectively.}
	
	\label{Fig:ExperimentC_2}
	
\end{figure}

\begin{figure}[t]
	\centering
	\subfigure[]
	{\label{Fig:ExperimentC_3a}
		\includegraphics[scale=0.34]{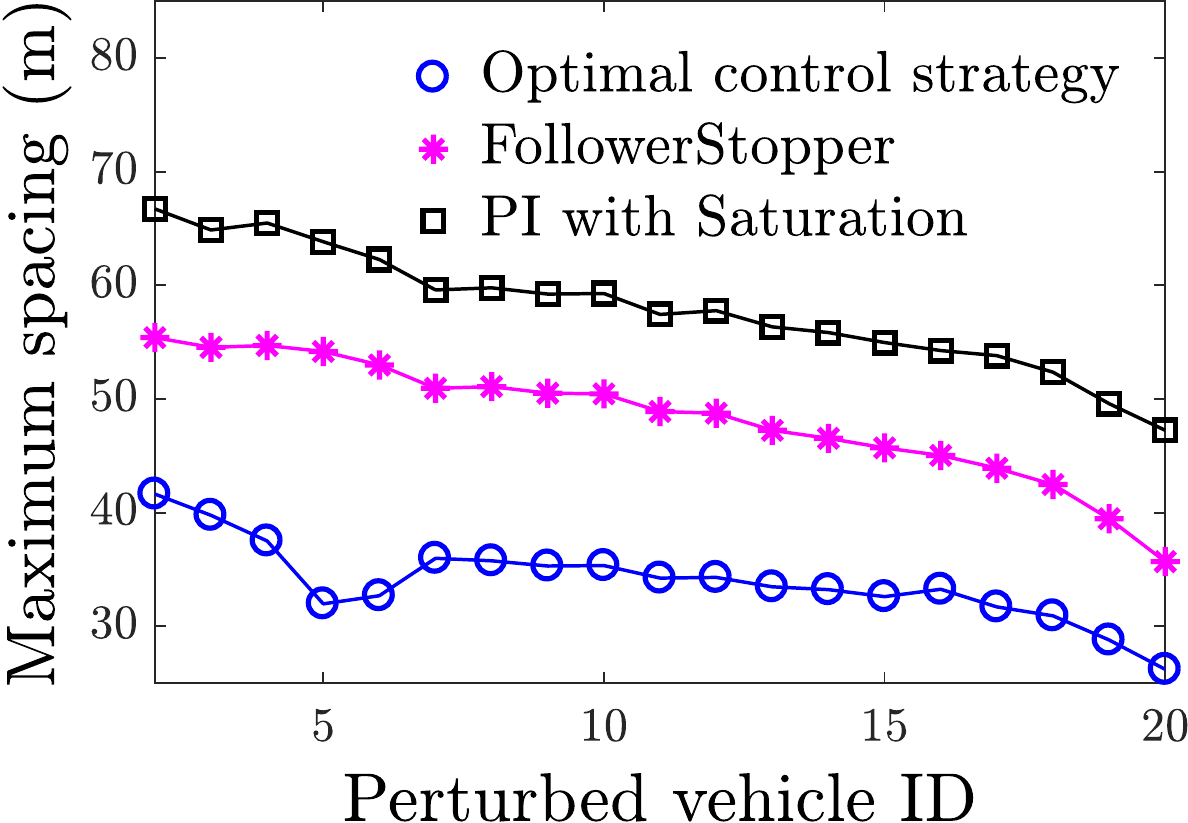}}
	\hspace{1mm}
	\subfigure[]
	{\label{Fig:ExperimentC_3d}
		\includegraphics[scale=0.34]{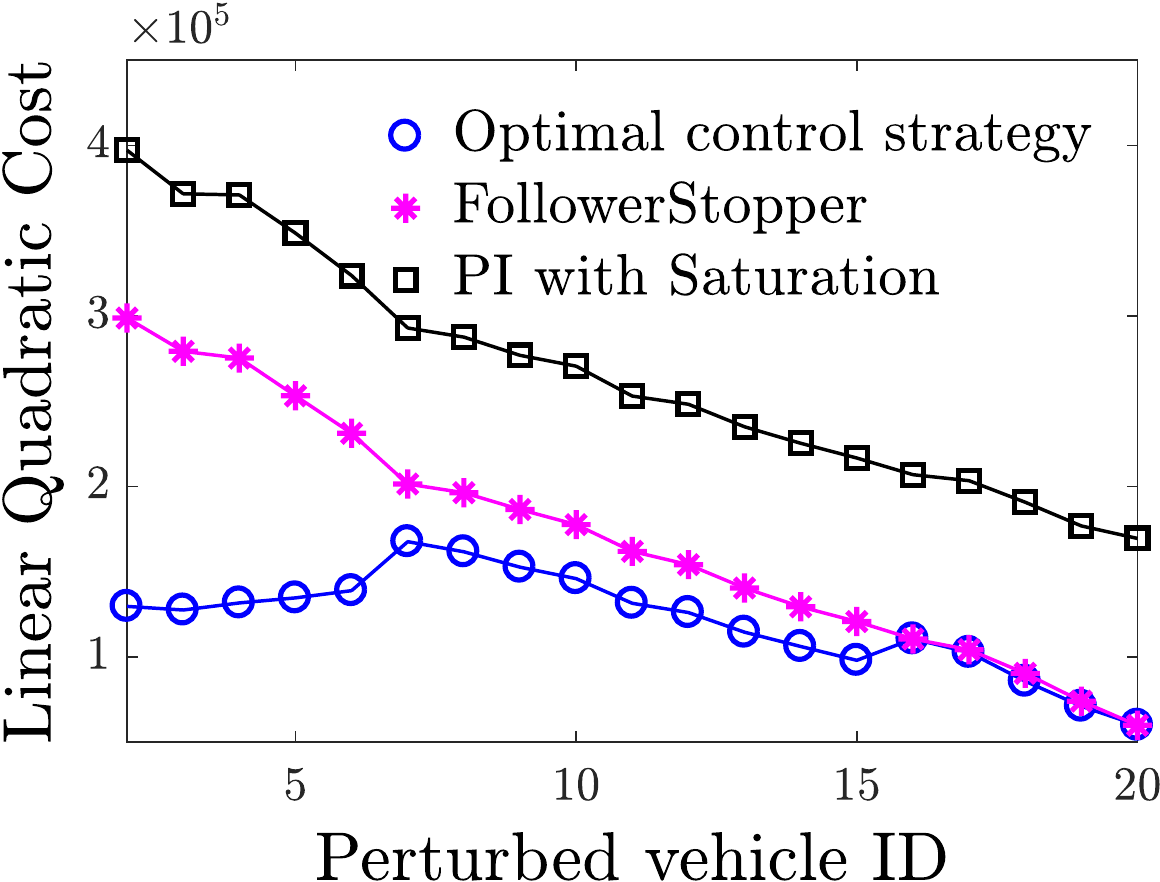}}
	\caption{Comparison of results at different positions of the perturbation (Experiment C). (a) The maximum spacing of the CAV during the overall process, \ie, $\max s_1(t)$. (b) The linear quadratic cost for the traffic system, defined as $\int_{t=0}^{\infty} x^{\tr}(t) Q x(t)+u^{\tr}(t) R u(t)$ with $Q$ and $R$ taking the same value as those in \eqref{Eq:Output}.}
	\label{Fig:ExperimentC_3}
\end{figure}

In our final experiment (Experiment C), we compare our proposed strategy with two existing heuristic ones: FollowerStopper and PI with Saturation \cite{stern2018dissipation}. Note that our experimental setup is not completely the same as that in \cite{stern2018dissipation}, and thus the performance of FollowerStopper and PI with Saturation in our experiments can be possibly improved after careful tuning. We assume that one vehicle is under a sudden perturbation, which may occur in the presence of infrastructure bottlenecks, \eg, intersections or merging lanes. The traffic flow has an initial velocity of $15\,\mathrm{m/s}$. At $t=20\,\mathrm{s}$, one HDV brakes at $-3\,\mathrm{m/s^2}$ for $3\,\mathrm{s}$. It is noticeably observed that when all the vehicles are under human control, the perturbed vehicle's action results in a traffic wave travelling upstream, which persists and does not vanish (Fig.~\ref{Fig:ExperimentC_1a}). If one CAV with the proposed method is included in the traffic system, it can prevent the propagation of the traffic wave and dampen the perturbation within a short time (Fig.~\ref{Fig:ExperimentC_1b}-\ref{Fig:ExperimentC_1d}).

The FollowerStopper and PI with Saturation in \cite{stern2018dissipation} use the command velocity $v_\mathrm{cmd}$ as the control input, and we add a proportional controller, \ie, $u(t)=k_p (v_\mathrm{cmd} (t)-v(t)$, with $k_p=0.6$, to serve as a lower controller. The comparison results are demonstrated in Fig.~\ref{Fig:ExperimentC_2} and Fig.~\ref{Fig:ExperimentC_3}. When all the vehicles are under human control, the spacing of each vehicle begins to fluctuate when the perturbation is introduced (Fig.~\ref{Fig:ExperimentC_2a}). If the CAV is using either of the three strategies, it is able to dampen the traffic wave and move the traffic state back to equilibrium (Fig.~\ref{Fig:ExperimentC_2b}-\ref{Fig:ExperimentC_2d}). However, the CAV leaves a very large spacing from the preceding vehicle for some time, when using FollowerStopper or PI with Saturation, with a maximum spacing larger than $50\,\mathrm{m}$ (Fig.~\ref{Fig:ExperimentC_2b}-\ref{Fig:ExperimentC_2c}). This gap may easily induce vehicles from adjacent lanes to cut in. In contrast, our optimal control strategy keeps the spacing within a moderate range during the whole transient process (Fig.~\ref{Fig:ExperimentC_2d}). Fig.~\ref{Fig:ExperimentC_3} illustrates two specific metrics with respect to the position of the perturbation. It is evident to see that our method witnesses an apparent reduction in terms of the maximum spacing of the CAV, which means our optimal controller overcomes the common shortage of previous methods that they tend to leave a large gap from the preceding vehicle. Besides, in terms of the overall cost during the system evolution, our method achieves the best performance, which is expected since our framework directly considers a system-level performance and attempts to obtain a corresponding optimal controller.

\section{Conclusion}
\label{Sec:6}

In this paper, we have proved that the ring-road mixed traffic system with one single CAV and multiple heterogeneous HDVs is stabilizable under a very mild condition. This provides a deeper insight towards the potential of traffic control via CAVs. In addition, we have established a theoretical framework to design a system-level optimal controller for mixed traffic systems, which utilizes the technique of structured optimal control to address the limit of the CAV's communication ability. Moreover, we have revealed that the CAV's desired spacing needs to be designed carefully since all the HDVs are not directly under control. Note that we focus on a ring road setup in this paper to show the stabilizability of mixed traffic systems and how to control traffic flow via CAVs. After modifying the system model \eqref{Eq:SystemModel}, the formulation of structured optimal control can be applied to a straight road scenario. However, the controllability results may be different since there is no constraint of ring-road structure and the corresponding uncontrollable mode is not relevant.

One future direction is to investigate the influence of model mismatch and drivers' reaction time on the mixed traffic system. These issues have been partially discussed in the estimation of HDV dynamics model~\cite{jin2018connected}, design of robust control strategies~\cite{jin2018experimental}, and stability analysis of mixed traffic systems~\cite{xie2018Heterogeneous}. Moreover, we consider a second-order model \eqref{Eq:LinearCAVModel} for the CAV, since the existing car-following models, \eg, OVM \cite{bando1995dynamical} and IDM \cite{Treiber2000Congested}, are mostly in a second-order form. It has been shown that the vehicle longitudinal dynamics play a significant role in cooperative control of multiple CAVs \cite{li2017dynamical}. Therefore, incorporating more practical dynamics model for the CAV is worth further investigation in the research of mixed traffic systems. Finally, considering that more than one CAVs may coexist in the traffic flow, another interesting topic is to design cooperative strategies for multiple CAVs to smooth traffic flow.


%

\appendices
\section{Proof of \eqref{Eq:Algebraa}}
\label{Sec:AppendixA}

Consider the solution of $\hat{A}^2 p=0$. Express $\hat{A}^2$ in a compact form as follows
$$
\hat{A}^2=\begin{bmatrix} F_{11} & 0 & \ldots &\ldots & F_{13} & F_{12} \\
F_{22} & F_{21} & 0 & \ldots & \ldots & F_{23} \\
F_{33} & F_{32} & F_{31} & 0 & \ldots & 0\\
\vdots & \ddots & \ddots & \ddots & \ddots & \vdots\\
0 & \ldots & \ldots & F_{n3} & F_{n2} & F_{n1}\end{bmatrix},
$$
where
$$
\begin{aligned}
F_{i1}&=A_{i1}^2,\\
F_{i2}&=A_{i2}A_{(i-1)1}+A_{i1}A_{i2},\\
F_{i3}&=A_{i2}A_{(i-1)2}.
\end{aligned}
$$
Similarly, $\hat{A}^2 p=0$ is equivalent to $F_{i1} p_i+F_{i2 }p_{i-1}+F_{i3} p_{i-1}=0$. For simplicity, we denote $d_i=\alpha_{i1} p_{i1}-\alpha_{i2} p_{i2}+\alpha_{i3} p_{(i-1)2}$, $i=1,\ldots,n$. Then the expanded expression of $F_{i1} p_i+F_{i2} p_{i-1}+F_{i3} p_{i-1}=0$ can be written as
\begin{subequations}\label{Eq:AhatpEqual0}
\begin{numcases}{}
d_{i-1}=d_i, \label{Eq:AhatpEqual0a}\\
-\alpha_{i2}d_i+\alpha_{i3}d_{i-1}=\alpha_{i1}\left(p_{i2}-p_{(i-1)2}\right),\label{Eq:AhatpEqual0b}
\end{numcases}
\end{subequations}
for $i=1,\ldots,n$. From \eqref{Eq:AhatpEqual0a} we know that $d_i=d$, $i=1,\ldots,n$ with $d$ denoting a constant. Substituting it into \eqref{Eq:AhatpEqual0b}, we have $p_{i2}-p_{(i-1)2}=\frac{\alpha_{i3}-\alpha_{i2}}{\alpha_{i1}}  d$, and then,
\begin{equation}
\label{Eq:Solutionofd}
d\sum_{i=1}^{n}\frac{\alpha_{i3}-\alpha_{i2}}{\alpha_{i1}} =\sum_{i=1}^{n}\left(p_{i2}-p_{(i-1)2}\right)=0.
\end{equation}
Indeed we know $\alpha_{i1}>0, \alpha_{i2}>\alpha_{i3}>0$ due to the real driving behavior of human drivers \cite{cui2017stabilizing}, and hence $ \sum_{i=1}^{n}\frac{\alpha_{i3}-\alpha_{i2}}{\alpha_{i1}}<0$. Therefore, $d_i=d=\alpha_{i1} p_{i1}-\alpha_{i2} p_{i2}+\alpha_{i3} p_{(i-1)2}=0$ is obtained from \eqref{Eq:Solutionofd}, and also, we have $p_{i2}-p_{(i-1)2}=0$ from \eqref{Eq:AhatpEqual0b}. These two results show that \eqref{Eq:AhatpEqual0} is equivalent to \eqref{Eq:ApzeroSolution}, indicating that the solution of $\hat{A}^2 p=0$ should also satisfy the form in \eqref{Eq:Solutionofp}. Thus, the dimension of the solution space of $\hat{A}p=0$ is one, \ie, $\dim \, \ker (\hat{A}-0\cdot I)^2=1 $. This completes the proof of \eqref{Eq:Algebraa}.

\section{Proof of \eqref{Eq:Eigenvalue}}
\label{Sec:AppendixB}
Consider a non-zero eigenvalue $\lambda$ of $\hat{A}$ and its left eigenvector $\rho$ with a same expression as that in \eqref{Eq:FormofRho}. Since $\rho^\tr (\lambda I-\hat{A} )=0$, we have ($i=1,\ldots,n$)
\begin{align}
	&\rho_{i1}\lambda - \rho_{i2}\alpha_{i1}=0, \label{Eq:RhoValuea}\\
	&\rho_{i1}+\rho_{i2}(\lambda+\alpha_{i2})-\rho_{(i+1)1}-\rho_{(i+1)2}\alpha_{(i+1)3}=0. \label{Eq:RhoValueb}
\end{align}
From \eqref{Eq:RhoValuea}, we have $\rho_{i2}=\rho_{i1}\cdot \frac{\lambda}{\alpha_{i1}}$. Substituting it into \eqref{Eq:RhoValueb} yields
\begin{equation}\label{Eq:EigenvalueEigenvectorEquation}
\left(1+\frac{\lambda^2}{\alpha_{i1}}+\frac{\alpha_{i2}}{\alpha_{i1}}\lambda\right)\rho_{i1}=\left(\frac{\alpha_{(i+1)3}}{\alpha_{(i+1)2}}\lambda+1\right)\rho_{(i+1)1},
\end{equation}
for $i=1,\ldots,n$. Assume there exists $k\in\{1,2,\ldots,n\}$, such that $\lambda^2+\alpha_{k2} \lambda+\alpha_{k1}=0$. In this case, the inequality, $\frac{\alpha_{(i+1)3}}{\alpha_{(i+1)2}}+1\neq 0$, $\forall i\in \{1,2,\ldots,n\}$, must hold; otherwise condition \eqref{Eq:ControllabilityTheorem} will be contradicted. Letting $i=k$ and substituting $\lambda^2+\alpha_{k2} \lambda+\alpha_{k1}=0$ into \eqref{Eq:EigenvalueEigenvectorEquation}, we have $\rho_{(k+1)1}=0$. Next, still applying \eqref{Eq:EigenvalueEigenvectorEquation} and letting $i=k+1$, we then have $\rho_{(k+2)1}=0$. After recursive iteration using \eqref{Eq:EigenvalueEigenvectorEquation}, we can obtain $\rho_{i1}=\rho_{i2}=0$, $i\in\{1,2,\ldots,n\}$, leading to $\rho=0$, which contradicts the requirement that $\rho \neq 0$. Accordingly, the assumption does not hold and therefore, $\lambda^2+\alpha_{i2} \lambda+\alpha_{i1}\neq0$, $\forall i\in\{1,2,\ldots,n\}$. Similarly, it can also be proved that $\alpha_{i3} \lambda +\alpha_{i2} \neq 0$, $\forall i\in\{1,2,\ldots,n\}$, based on \eqref{Eq:EigenvalueEigenvectorEquation}.



\ifCLASSOPTIONcaptionsoff
  \newpage
\fi



%

\bibliographystyle{IEEEtran}
\bibliography{IEEEabrv,mybibfile}

%

\end{document}

%% file: defs.tex
\def\B{{\bf B}}

\def\R{{\bf R}}

\def\Z{{\bf Z}}

\def\0{{\bf 0}}
\def\1{{\bf 1}}

\def\etal{{\em et al.}}
\def\eg{{\em e.g.}}
\def\ie{{\em i.e.}}

\def\tr{\mathrm{T}}
\def\dim{\mathrm{dim}}
\def\ker{\mathrm{ker}}
\def\rank{\mathrm{rank}}
\def\etal{{\em et al.\/}\,}
